\documentclass[sn-mathphys,Numbered]{sn-jnl}

\usepackage{graphicx}%
\usepackage{multirow}%
\usepackage{amsmath,amssymb,amsfonts}%
\usepackage{amsthm}%
\usepackage{amsmath}
\usepackage{mathrsfs}%
\usepackage[title]{appendix}%
\usepackage{xcolor}%
\usepackage{textcomp}%
\usepackage{manyfoot}%
\usepackage{booktabs}%
\usepackage{algorithm}%
\usepackage{algorithmicx}%
\usepackage{algpseudocode}%
\usepackage{listings}%
\usepackage{indentfirst}
\usepackage{appendix}

\usepackage{subcaption} 
\usepackage{bm}

\newtheorem{theorem}{Theorem}
\newtheorem{lemma}[theorem]{Lemma}

\theoremstyle{remark}
\newtheorem{remark}{Remark}%
\newtheorem{assumption}{Assumption}

\begin{document}

\title[Article Title]{An Iteratively Decoupled Algorithm for Multiple-Network Poroelastic Model with Applications in Brain Edema Simulations} 


\author*[1]{\fnm{Mingchao} \sur{Cai}}\email{cmchao2005@gmail.com}

\author[2]{\fnm{Meng} \sur{Lei}}\email{210902070@njnu.edu.cn}

\author[3]{\fnm{Jingzhi} \sur{Li}}\email{li.jz@sustech.edu.cn}

\author[3]{\fnm{Jiaao} \sur{Sun}}\email{sunja@sustech.edu.cn}

\author[2]{\fnm{Feng} \sur{Wang}}\email{fwang@njnu.edu.cn}

\affil*[1]{\orgdiv{Department of Mathematics}, \orgname{Morgan State University}, \orgaddress{\city{Baltimore}, \postcode{21251}, \state{Maryland}, \country{USA}}}

\affil[2]{\orgdiv{Key Laboratory of NSLSCS, Ministry of Education; Jiangsu International Joint Laboratory of BDMCA; School of Mathematical Sciences}, \orgname{Nanjing Normal University}, \orgaddress{\city{Nanjing}, \postcode{210023}, \state{Jiangsu}, \country{P.R. China}}}

\affil[3]{\orgdiv{Department of Mathematics}, \orgname{Southern University of Science and Technology}, \orgaddress{\city{Shenzhen}, \postcode{21251}, \state{Guangdong}, \country{P. R. China}}}




\abstract{
In this work, we present an iteratively decoupled algorithm for solving the quasi-static multiple-network poroelastic model. 
Our approach employs a total-pressure-based formulation with solid displacement, total pressure, and network pressures as primary unknowns. This reformulation decomposes the original problem into a generalized Stokes problem and a parabolic problem, offering key advantages such as reduced elastic locking effects and simplified discretization. The algorithm guarantees unconditional convergence to the solution of the fully coupled system.  
Numerical experiments demonstrate the accuracy, efficiency, and robustness of the method with respect to physical parameters and discretization. We further apply the algorithm to simulate the brain edema process, showcasing its practical utility in biomechanical modeling. 
}



\maketitle

\section{Introduction}\label{sec1}

Let $\Omega \subset \mathbb{R}^d (d=2,3)$ be a bounded polygonal domain. The quasi-static multiple-network poroelastic problem \cite{hong2020parameter} is to find the displacement $\pmb{u}$ and the network pressures $\vec{p}$ such that
\begin{subequations} \label{strong_form12}
\begin{align}
\label{strong_form_eq1}
&-\operatorname{div} \sigma(\pmb{u})+ \nabla (\vec{\alpha}^{\mathsf{T}} \vec{p}) =\boldsymbol{f} \quad \text { in } \Omega \times[0, T], \\
\label{strong_form_eq2}
& \vec{\alpha} \operatorname{div} \dot{\pmb{u}} +S \dot{\vec{p}}+B\vec{p}-\operatorname{div}(K \nabla \vec{p})=\vec{g} \quad \text { in } \Omega \times[0, T].
\end{align}
\end{subequations}
Here, $\pmb{u}=(u_1(\pmb{x},t),u_2(\pmb{x},t))^{\mathsf{T}}$\,and $\vec{p}=(p_1(\pmb{x},t),p_2(\pmb{x},t),\cdots,p_N(\pmb{x},t))^{\mathsf{T}}$, with a given number of networks $N$.
The operators and parameters are defined as follows. The effective stress and the strain tensor are denoted by
    $\sigma(\pmb{u})=2\mu \epsilon(\pmb{u})+\lambda \operatorname{div} (\pmb{u})\pmb{I}$ and
    $\epsilon(\pmb{u})=\frac{1}{2}\left(\nabla \pmb{u}+(\nabla \pmb{u})^{\mathsf{T}}\right)$
 , respectively. The Lam\'{e} parameters $\lambda$ and $\mu$ are expressed in terms of the Young modulus $E$ and Poisson ratio $\nu \in [0,\frac{1}{2})$ by $
    \lambda=\frac{\nu E}{(1+\nu)(1-2\nu)}$ and
    $\mu=\frac{E}{2(1+\nu)}$, respectively.
Column vector $\vec{\alpha}=(\alpha_1,\alpha_2,\cdots,\alpha_N)^{\mathsf{T}}$, where $\alpha_i\in (0,1]$ is the  Biot-Willis coefficient. Matrix $S=\text{diag}(c_{1},c_{2},\cdots,c_{N})$, and $c_{i}\geq 0$ is the storage coefficient. Matrix $K=\text{diag}(K_1,K_2,\cdots,K_N)$, and $K_i$ is the hydraulic conductivity coefficient. 
The coefficient matrix $B$ satisfies $(B\vec{p})_i=\sum_{j=1, j \neq i}^N \beta_{i j}\left(p_i-p_j\right)$, where the non-negative network transfer coefficient $\beta_{ij}$ couple the network pressures and $\beta_{ij}=\beta_{ji}, 1\leq i,j\leq N, j\neq i$. Further, $\pmb{f}=(f_1,f_2)^{\mathsf{T}}$, where $f_i$ represents the body force and $\vec{g}=(g_1,g_2,\cdots,g_N)^{\mathsf{T}}$, where $g_i$ represents the source/sink in the $i$th network.

The system \eqref{strong_form12} is well-posed with proper boundary and initial conditions \cite{hong2020parameter, lee2017parameter}. In this paper, we consider the mixed boundary conditions, assuming $\partial \Omega=\Gamma_{\pmb{u}, D} \cup \Gamma_{\pmb{u}, N}=\Gamma_{\vec{p}, D} \cup \Gamma_{\vec{p}, N},|\Gamma_{\pmb{u}, D} \cap \Gamma_{\pmb{u},N}|=0,|\Gamma_{\vec{p},D} \cap \Gamma_{\vec{p}, N}|=0,|\Gamma_{\pmb{u},D}|>0,|\Gamma_{\vec{p}, D}|>0$. More clearly, the following boundary and initial conditions \cite{lee2019mixed} are imposed.
\begin{align*}
    \pmb{u}&=\pmb{0}\quad  \text{on} \, \Gamma_{\pmb{u},D},\\
    \left(\sigma(\pmb{u})- (\vec{\alpha}^{\mathsf{T}} \vec{p}) \pmb{I}\right)\pmb{n}&=\pmb{h}\quad  \text{on} \, \Gamma_{\pmb{u},N},\\
    \vec{p} &=\vec{0} \quad   \text{on} \, \Gamma_{\vec{p},D},\\
    \left( K \nabla \vec{p} \right) \pmb{n}&=\vec{l} \quad  \text{on} \, \Gamma_{\vec{p},N},\\
    \pmb{u}(0)&=\pmb{u}_0 \quad  \text{in} \, \Omega, \\
    \vec{p}(0)&=\vec{p}_{0} \quad  \text{in} \, \Omega,
\end{align*}
where $\pmb{u}_0$ and $\vec{p}_0$ are assumed to satisfy the \eqref{strong_form_eq1}.

Multiple-network poroelasticity, introduced in geomechanics \cite{bai1993multiporosity}, models mechanical deformation and fluid flow in multiple porous materials. It extends Biot's model, which is widely used in geomechanics, geophysics, and biology. For instance, when considering the case with $N=2$, one can obtain the Biot-Barenblatt model \cite{nordbotten2010posteriori,showalter2002single}, which states consolidation processes within a fluid-saturated double-diffusion model of fractured rock.
In another context, Tully and Ventikos investigated a scenario involving four distinct networks ($N=4$) as a macroscopic model to describe the dynamics of fluid flows within brain tissue \cite{tully2011cerebral}.
More recently, a multiple-network poroelasticity model has also been employed to gain a deeper understanding of the impact of biomechanical risk factors associated with the early stages of Alzheimer's disease \cite{guo2018subject}.

The discretization of Biot's equations is widely recognized as a challenging task due to the poroelastic locking phenomenon. Poroelastic locking is characterized by two main features \cite{lee2019mixed}: (1) underestimation of solid deformation if the material is close to being incompressible and (2) non-physical pressure oscillations.  Many existing studies, such as \cite{both2017robust, mikelic2014numerical}, adopt solid displacement and fluid pressure as the primary variables within Biot's consolidation model. However, it has been observed that the elastic locking phenomenon occurs in a model based on the two-field formulation \cite{cai2015comparisons, lee2017parameter}.
To address these challenges, researchers have employed various approaches, including the DG method \cite{phillips2009overcoming}, stabilization techniques \cite{rodrigo2016stability, rodrigo2018new, kim2011stability, mikelic2013convergence}, and various three-field or four-field reformulations \cite{hong2017parameter, hu2017nonconforming, ju2020parameter, lee2016robust, oyarzua2016locking, qi2021finite, wang2022mixed, yi2017iteratively}.
In spite of the extensive research dedicated to exploring Biot's consolidation model, the
 multiple-network poroelasticity problem has received much less attention. Hong et al. introduced parameter-robust preconditioners \cite{hong2019conservative}, proposed a fixed-stress split method \cite{hong2020parameter} and presented an augmented Lagrangian Uzawa algorithm \cite{hong2020parameterUzawa} for multiple-network flux-based
poroelasticity equations. Lee et al. analyzed a mixed finite element method \cite{lee2019mixed} and proposed two partitioned numerical algorithms  \cite{lee2019unconditionally} for the multiple-network poroelasticity problem.

By introducing an intermediate variable $\xi=\vec{\alpha}^{\mathsf{T}} \vec{p}-\lambda \operatorname{div} \pmb{u}$, we studied a formulation which has the
solid displacement, the total pressure, and the network pressures as the primary unknowns. This formulation interprets the original problem as a combination of the generalized Stokes problem and the parabolic problem. The solid displacement and the total pressure are approximated by the classical stable Stokes finite elements, while the network pressures are discretized by the Lagrange finite elements. Based on this formulation, we propose an iteratively decoupled algorithm for solving the multiple-network poroelasticity problem \eqref{strong_form12}. Unlike conventional two-field formulation approaches, our algorithm offers two key advantages: (1) it overcomes the difficulties caused by Poisson locking and requires no stabilization parameters, and (2) it guarantees unconditional convergence to the coupled solution. 
The algorithm remains robust even in limiting cases, such as when storage coefficients satisfy $c_{i} = 0$ for $i=1,2,\cdots, N$. Numerical experiments confirm that the method maintains optimal convergence order across varying physical parameters. 
To further validate the algorithm's practical efficacy, we apply it to brain edema simulations. By comparing the results with those obtained from the coupled algorithm, we illustrate that the decoupled algorithm provides the same order of accuracy while requiring significantly fewer computational resources.

The structure of this article is as follows. In Section 2, we present the total-pressure-based formulation  
and the corresponding variational problems for the multiple-network poroelasticity equations.
The iteratively decoupled algorithm is given in Section 3. In Section 4, we analyze the convergence of the iteratively decoupled algorithm. 
Numerical experiments and brain edema simulations are presented in Section 5. Conclusions are drawn in Section 6.

\section{A total-pressure-based formulation}



We employ the notation $L^2(\Omega)$ and $L^2(\partial \Omega)$ to represent real-valued functions that are square-integrable over the domains $\Omega$ and $\partial \Omega$, respectively. The inner product for the functions of $L^2(\Omega)$ is denoted by $(\cdot, \cdot)$ with its corresponding induced norm represented as $|\!|\cdot|\!|_{L^2(\Omega)}$. Similarly, for vector-valued or matrix-valued functions, where each component belongs to $L^2(\Omega)$, we utilize the same notation $(\cdot, \cdot)$ and $|\!|\cdot|\!|_{L^2(\Omega)}$ to express the inner product and norm. In the case of vector-valued functions, where each component belongs to $L^2(\partial \Omega),$ we also employ the notation $\langle \cdot, \cdot \rangle$ to denote the inner product.
For a non-negative integer $m$ and $1 < p < +\infty$, we denote by 
$W^{m,p}(\Omega)=\{ u |D^{\alpha}u \in L^p (\Omega), 0 \leq |\alpha| \leq m, |\!|u|\!|_{W^{m,p}(\Omega)} < +\infty\}$ the Sobolev spaces; When $p=2$, we  use $H^m(\Omega)$ to denote
$W^{m,2}(\Omega)$. 
In cases involving vector-valued functions where each component belongs to $H^m(\Omega)$, we consistently employ the notation $|\!|\cdot|\!|_{H^m(\Omega)}$ to represent the norm.
Furthermore, when $m\geq 1$, we use $H^m_{0,\Gamma}(\Omega)$ to denote the subspace of $H^m (\Omega)$ comprising functions with a vanishing trace on $\Gamma \subset \partial \Omega$.
For ease of presentation, we assume that $\Omega$ is a two-dimensional domain in this work.


In the following, we introduce an intermediate variable to reformulate the problem \eqref{strong_form12} and present the associated variational problem \cite{lee2019mixed}. Specifically, we define the so-called "total pressure" as $\xi=\vec{\alpha}^{\mathsf{T}}\vec{p}-\lambda \operatorname{div} \pmb{u}$. With this, \eqref{strong_form12} can be expressed as
\begin{subequations} \label{three_field_form}
\begin{align}
\label{three_field_form1}
&-2 \mu \operatorname{div}\left( \varepsilon(\boldsymbol{u})\right) +\nabla \xi=\boldsymbol{f}, \\
\label{three_field_form2}
&-\operatorname{div} \boldsymbol{u}-\frac{1}{\lambda} \xi+\frac{1}{\lambda} \vec{\alpha}^{\mathsf{T}}\vec{p} =0, \\
\label{three_field_form3}
&\left(S+\frac{1}{\lambda} \vec{\alpha}\vec{\alpha}^{\mathsf{T}}  \right)\dot{\vec{p}}-\frac{1}{\lambda}\vec{\alpha}
\dot{\xi} -\operatorname{div}(K\nabla \vec{p}) +B\vec{p}=\vec{g}.
\end{align}
\end{subequations}
The initial and boundary conditions we provided earlier can still be applied to the problem \eqref{three_field_form}.

In order to describe the variational problem for the total-pressure-based formulation  \eqref{three_field_form}, we will utilize the following functional spaces.
\begin{align*}
   \pmb{V}&:=\{\pmb{v}\in [H^1 (\Omega)]^2;\pmb{v}|_{ \Gamma_{\pmb{u},D}} =\pmb{0}\},\\ \notag
    W&:=L^2 (\Omega),\\  \notag
    M&:=\{\vec{q}\in [H^1(\Omega)]^N; \vec{q} |_{\Gamma_{\vec{p},D}}=\vec{0}\}.
\end{align*}
The properties of these spaces are as follows.
The first Korn inequality \cite{gu2023iterative} holds on $\pmb{V}$, that is, there exists a constant $C_K=C_K (\Omega,\Gamma_{\pmb{u},D})>0$ such that
$$|\!|\pmb{u}|\!|_{H^1(\Omega)} \leq C_K |\!|\epsilon(\pmb{u})|\!|_{L^2(\Omega)} \quad \forall \, \pmb{u}\in \pmb{V}.$$
Furthermore, the following inf-sup condition  \cite{gu2023iterative} holds. There exists a constant $\beta_0>0$ depending only on $\Omega$ and $\Gamma_{\pmb{u},D}$ such that
\begin{align} \label{inf-sup condition}
\sup_{\pmb{v}\in \pmb{V}} \frac{(\operatorname{div} \pmb{v},\eta)}{|\!|\pmb{v}|\!|_{H^1 (\Omega)}} \geq \beta_0 |\!|\eta|\!|_{L^2(\Omega)} \quad \forall \,\eta \in W.
\end{align}

\begin{assumption} \label{Assumption1}
    We assume that $\pmb{u}_0 \in \pmb{H}^1(\Omega)$, $\pmb{f}\in \pmb{L}^2(\Omega)$, $\pmb{h}\in \pmb{L}^2(\Gamma_{\pmb{u},N})$, $\vec{p}_0 \in [L^2(\Omega)]^N$, $\vec{g} \in [L^2(\Omega)]^N$  and $\vec{l}\in [L^2({\Gamma_{\vec{p},N}})]^N$. We also assume that $\mu>0$, $\lambda>0$, $K_i>0, i=1,\cdots,N$, $\beta_{i,j}>0, 1\leq i,j\leq N ,j\neq i$, $T>0$.
\end{assumption}

For simplicity, we will assume \textbf{Assumption} \ref{Assumption1} holds in the rest of our paper. For ease of presentation, we also assume that $\pmb{f}$, $\pmb{h}$, $\vec{g}$, $\vec{l}$ are independent of $t$. Given $T>0$, $(\pmb{u}, \xi, \vec{p})\in \pmb{V}\times W \times M$ with
\begin{align*}
   &\pmb{u}\in L^{\infty}(0,T;\pmb{V}), \, \xi\in L^{\infty}(0,T;W), \, \vec{p}\in L^{\infty}(0,T;M),\\
   & \dot{\vec{p}}\in L^2(0,T;M'), \, \dot{\xi}\in L^2(0,T;M'),
\end{align*}
is called a weak solution of problem \eqref{three_field_form1}-\eqref{three_field_form3}, if there holds
\begin{subequations} \label{weak_problem}
\begin{align}
\label{three_field_form_wp1}
&2 \mu (\epsilon(\pmb{u}), \epsilon(\pmb{v}))-(\xi, \operatorname{div} \pmb{v}) =(\pmb{f}, \pmb{v})+\langle \pmb{h},\pmb{v}\rangle_{\Gamma_{\pmb{u},N}},\\
\label{three_field_form_wp2}
&(\operatorname{div} \pmb{u}, \eta )+\frac{1}{\lambda}( \xi, \eta)-\frac{1}{\lambda}( \vec{\alpha}^{\mathsf{T}}\vec{p}, \eta)=0,\\
\label{three_field_form_wp3}
 \left(  \left(S+\frac{1}{\lambda} \vec{\alpha}\vec{\alpha}^{\mathsf{T}}  \right)\dot{\vec{p}},\vec{q} \right)& -\frac{1}{\lambda}\left( \vec{\alpha}\dot{\xi},\vec{q} \right)+\left(  K\nabla \vec{p},\nabla \vec{q}   \right)+\left( B\vec{p},\vec{q}\right)=\left( \vec{g},\vec{q}\right)+\left\langle \vec{l},\vec{q}\right\rangle_{\Gamma_{\vec{p},N}},
 \end{align}
\end{subequations}
for almost every $t\in [0,T]$.

\section{The algorithms}

In this section, building upon the previously described total-pressure reformulation, we will introduce both the coupled algorithm and the iteratively decoupled algorithm. Our discussion will commence with an explanation of the mixed finite element spaces used to discretize the variational problem \eqref{weak_problem}. 

Let $\mathcal{T}_h$ denote a shape-regular and conforming triangulation of the domain $\Omega$. Define $h:= \max_{K \in \mathcal{T}_h} h_K$.
The solid displacement and total pressure are approximated by the lowest-order Taylor-Hood elements, i.e.,$(\pmb{P}_2, P_1)$ Lagrange finite elements. In contrast, the network pressures
are discretized by using the $P_1$ Lagrange finite elements \cite{taylor1973numerical}.
The finite element spaces on $\mathcal{T}_h$ are as follows.
\begin{align*}
&\pmb{V}_h:=\{\pmb{v}_h \in \pmb{V}\cap  [C^0(\overline{\Omega})]^2;\pmb{v}_h|_{\Gamma_{\pmb{u},D}} =\pmb{0},\pmb{v}_h |_K \in [P_2(K)]^2,\forall \,K \in T_h\},\\
&W_h :=\{\eta_h \in W \cap C^0(\overline{\Omega});\eta_h |_K \in P_1(K),\forall \, K \in T_h\},\\
&M_h:=\{\vec{q}_{h}\in  M \cap [C^0(\overline{\Omega})]^N;\vec{q}_h|_{ \Gamma_{\vec{p},D}}=\vec{0},\vec{q}_h |_K\in [P_1(K)]^N,\forall \, K\in T_h\}.
\end{align*}
We note that $\pmb{V}_h \times W_h$ is a stable Stokes pair, i.e., the following discrete inf-sup condition holds. There exists a constant $\beta_0^* >0$, independent of $h$, such that
\begin{align} \label{discrete inf-sup condition}
    \sup_{\pmb{v}_h \in \pmb{V}_h} \frac{(\operatorname{div} \pmb{v}_h,\eta_h)}{|\!|\pmb{v}_h|\!|_{H^1(\Omega)}}\geq  \beta_0^* |\!|\eta_h|\!|_{L^2(\Omega)}\quad \forall \, \eta_h \in W_h.
\end{align}

\subsection{A coupled algorithm}
The time discretization is taken as an equidistant partition $0=t^0<t^1<\cdots <t^L=T$. For simplicity, we define $\pmb{u}^n=\pmb{u}(t^n),\xi^n=\xi(t^n),\vec{p}^n=\vec{p}(t^n)$.
Suppose that initial values $(\pmb{u}_h^0,\xi_h^0,\vec{p}_{h}^0)\in \pmb{V}_h \times W_h \times M_{h}$ are provided, we apply a backward Euler
scheme for the time discretization. The following algorithm is considered.

For all $n=1,2,\cdots,L$, given $(\pmb{u}_h^{n-1},\xi_h^{n-1},\vec{p}_{h}^{n-1})\in \pmb{V}_h \times W_h \times M_{h}$, \,find $(\pmb{u}_h^n,\xi_h^n,\vec{p}_{h}^n)\in \pmb{V}_h \times W_h \times M_{h}$, for all $(\pmb{v}_h,\eta_h,\vec{q}_{h})\in \pmb{V}_h \times W_h \times M_{h}$ such that
\begin{subequations} \label{coupled_alg}
\begin{align} \label{coupled_alg1}
2\mu (\epsilon(\pmb{u}_h^n),\epsilon(\pmb{v}_h))&-(\xi_h^n,\operatorname{div} \pmb{v}_h)=(\pmb{f},\pmb{v}_h)+\langle \pmb{h},\pmb{v}_h \rangle_{\Gamma_{\pmb{u},N}},\\
\label{coupled_alg2}
(\operatorname{div}\pmb{u}_h^n,\eta_h)&+\frac{1}{\lambda}(\xi_h^n,\eta_h)-\frac{1}{\lambda}(\vec{\alpha}^{\mathsf{T}} \vec{p}_h^n,\eta_h)=0,\\
\label{coupled_alg3}
  \left(  \left(S+\frac{1}{\lambda} \vec{\alpha}\vec{\alpha}^{\mathsf{T}}  \right)\frac{\vec{p}_h^n-\vec{p}_h^{n-1}}{\Delta t},\vec{q}_h \right) &-\frac{1}{\lambda}\left( \vec{\alpha}\frac{\xi_h^n-\xi_h^{n-1}}{\Delta t},\vec{q}_h \right)+\left(  K\nabla \vec{p}_h^n,\nabla \vec{q}_h   \right)+\left( B\vec{p}_h^n,\vec{q}_h\right) \\ \notag
=\left( \vec{g},\vec{q}_h\right)+\left\langle \vec{l},\vec{q}_h\right\rangle_{\Gamma_{\vec{p},N}}.
\end{align}
\end{subequations}

\subsection{An iteratively decoupled algorithm} We propose an iteratively decoupled algorithm where, at each time step $t^n$, the reaction-diffusion equation for $\vec{p}$ and the generalized Stokes equations for $\pmb{u}$ and $\xi$ are solved alternately, using the previous iteration's solution to decouple the computations.

Let us define a sequence $\{(\pmb{u}_h^{n,k},\xi_h^{n,k},\vec{p}_h^{n,k})\}$ with $k \in \mathbb{N}^*$ being the iteration index. After the initialization, i.e., $\pmb{u}_h^{n,0}=\pmb{u}_h^{n-1},\xi_h^{n,0}=\xi_h^{n-1},\vec{p}_{h}^{n,0}=\vec{p}_{h}^{n-1}$, each iteration consists of  the following two steps.
For a fixed time-step index $n$, the $k$-th iteration can be expressed as follows.\\
 \textbf{Step\,1}. Given $\xi_h^{n,k-1}\in W_h$, find $\vec{p}_{h}^{n,k}\in M_{h}$, for all $\vec{q}_h \in M_h$ such that
\begin{subequations}
\begin{align} \label{decoupled_alg_step1}
  &\left(  \left(S+\frac{1}{\lambda} \vec{\alpha}\vec{\alpha}^{\mathsf{T}}  \right)\vec{p}_h^{n,k},\vec{q}_h \right) +\Delta t \left(  K\nabla \vec{p}_h^{n,k},\nabla \vec{q}_h   \right)+\Delta t\left( B\vec{p}_h^{n,k},\vec{q}_h\right)\\ \notag
  &=\left(  \left(S+\frac{1}{\lambda} \vec{\alpha}\vec{\alpha}^{\mathsf{T}}  \right)\vec{p}_h^{n-1},\vec{q}_h \right)
  + \frac{1}{\lambda}\left(\vec{\alpha}\left(\xi_h^{n,k-1}-\xi_h^{n-1}\right),\vec{q}_h \right)    +\Delta t\left( \vec{g},\vec{q}_h\right)+\Delta t\left\langle \vec{l},\vec{q}_h\right\rangle_{\Gamma_{\vec{p},N}}.
\end{align}\\
\textbf{Step\,2}. Given $\vec{p}_{h}^{n,k}\in M_{h}$, find $(\pmb{u}_h^{n,k},\xi_h^{n,k})\in \pmb{V}_h\times W_h$, for all $(\pmb{v}_h,\eta_h)\in \pmb{V}_h \times W_h$ such that
\begin{align} \label{decoupled_alg_step2}
   & 2\mu (\epsilon(\pmb{u}_h^{n,k}),\epsilon(\pmb{v}_h))-(\xi_h^{n,k},\operatorname{div}\pmb{v}_h)=(\pmb{f},\pmb{v}_h)+\langle \pmb{h},\pmb{v}_h\rangle_{\Gamma_{\pmb{u},N}},\\
 \label{decoupled_alg_step3}
  & (\operatorname{div} \pmb{u}_h^{n,k},\eta_h)+\frac{1}{\lambda}(\xi_h^{n,k},\eta_h)-\frac{1}{\lambda}(\vec{\alpha}^{\mathsf{T}} \vec{p}_h^{n,k},\eta_h)=0.
\end{align}
\end{subequations}

\section{Convergence analysis}

For the error analysis of the coupled algorithm, we direct the readers to \textbf{Appendix A}. In this section, our objective is to establish that as $k$ approaches positive infinity, the sequence ${(\pmb{u}_h^{n,k},\xi_h^{n,k},\vec{p}_h^{n,k})}$ generated through the iteratively decoupled algorithm converges towards the solution ${(\pmb{u}_h^n,\xi_h^n,\vec{p}_h^n)}$ obtained via the coupled algorithm. To substantiate these convergence results, we present the following lemma \cite{oyarzua2016locking, lee2017parameter, gu2023iterative}.



\begin{lemma}
    For all $\pmb{u}_h \in \pmb{V}_h$, \, the following inequality holds
 \begin{align} \label{Korn_conclusion}
|\!| \operatorname{div}\pmb{u}_h|\!|_{L^2(\Omega)}\leq \sqrt{d}|\!|\epsilon(\pmb{u}_h)|\!|_{L^2(\Omega)}.
\end{align}
\end{lemma}

For a complete theoretical analysis of the decoupled iterative algorithm, 
we firstly assume that $\min_{1\leq i \leq N} c_{i}>0$, and we have the following theorem. 

\begin{theorem}
    Let ${(\pmb{u}_h^n,\xi_h^n,\vec{p}_h^n)}$ and ${(\pmb{u}_h^{n,k},\xi_h^{n,k},\vec{p}_h^{n,k})}$ be the solutions of problem \eqref{coupled_alg1}-\eqref{coupled_alg3} and problem \eqref{decoupled_alg_step1}-\eqref{decoupled_alg_step3}, respectively. Let $e_{\pmb{u}}^k=\pmb{u}_h^{n,k}-\pmb{u}_h^n,e_{\xi}^k=\xi_h^{n,k}-\xi_h^n,e_{\vec{p}}^k=\vec{p}_{h}^{n,k}-\vec{p}_{h}^n$ denote the errors between the iterative solution in the $k$-th step and the solution of the coupled algorithm. Then, for all $k\geq 1$,  it holds that

\begin{align*} 
    |\!|e_{\xi}^k|\!|_{L^2(\Omega)} \leq C^* |\!|e_{\xi}^{k-1}|\!|_{L^2(\Omega)},
\end{align*}
where $C^*=\frac{\frac{|\!|\vec{\alpha}|\!|^2}{\lambda}}{\delta+\frac{|\!|\vec{\alpha}|\!|^2}{\lambda}}$ is a positive constant less than  1, \,$|\!|\vec{\alpha}|\!|=\left(\sum_{i=1}^N \alpha_i^2\right)^{\frac{1}{2}}$,\,$\delta=\min_{1\leq i\leq N} c_{i}$.
Moreover, it holds that
\begin{align*}
 \lim _{k\rightarrow +\infty} |\!|\nabla e_{\vec{p}}^k|\!|_{L^2(\Omega)}&=
     \lim _{k\rightarrow +\infty} |\!|e_{\vec{p}}^k|\!|_{L^2(\Omega)}=
    \lim _{k\rightarrow +\infty} |\!|e_{\pmb{u}}^k|\!|_{H^1(\Omega)}=
     \lim _{k\rightarrow +\infty} |\!|e_{\pmb{u}}^k|\!|_{L^2(\Omega)}=0. \notag
\end{align*}

\end{theorem}

\noindent{\bf Proof}. Subtracting \eqref{coupled_alg3}, \eqref{coupled_alg1} and \eqref{coupled_alg2} from \eqref{decoupled_alg_step1}, \eqref{decoupled_alg_step2} and \eqref{decoupled_alg_step3}, and taking $\pmb{v}_h =e_{\pmb{u}}^k,\eta_h=e_{\xi}^k,\vec{q}_{h}=e_{\vec{p}}^k$, respectively, we have
\begin{align}
\label{decoupled_alg_ana1}
  &  2\mu (\epsilon(e_{\pmb{u}}^k),\epsilon(e_{\pmb{u}}^k))-(e_{\xi}^k,\operatorname{div} e_{\pmb{u}}^k)=0,\\
\label{decoupled_alg_ana2}
& (\operatorname{div}e_{\pmb{u}}^k,e_{\xi}^k)+\frac{1}{\lambda}(e_{\xi}^k,e_{\xi}^k)=\frac{1}{\lambda}(\vec{\alpha}^{\mathsf{T}} e_{\vec{p}}^k,e_{\xi}^k),\\
\label{decoupled_alg_ana3}
\left( \left( S+\frac{1}{\lambda} \vec{\alpha}\vec{\alpha}^{\mathsf{T}} \right)e_{\vec{p}}^k,e_{\vec{p}}^k\right)&+\Delta t \left(  K\nabla e_{\vec{p}}^k,\nabla e_{\vec{p}}^k \right)+\Delta t \left(  Be_{\vec{p}}^k, e_{\vec{p}}^k \right)=\frac{1}{\lambda}\left(  \vec{\alpha}e_{\xi}^{k-1},e_{\vec{p}}^k \right).
\end{align}
Here,
\begin{align} \label{identity_1}
\left(B\vec{q},\vec{q}\right)= \sum_{1\leq i<j\leq N}\beta_{i,j} |\!|q_{i}-q_{j}|\!|_{L^2(\Omega)}^2 \quad \forall \, \vec{q}\in M.
\end{align}
From this expression and the definition of $K$, we know that $B$ is positive semi-definite and $K$ is positive definite. Combining with \eqref{decoupled_alg_ana3}, one can obtain
\begin{align} \label{decoupled_alg_ana8}
     \left(  Se_{\vec{p}}^k, e_{\vec{p}}^k\right) +\frac{1}{\lambda} |\!|\vec{\alpha}^{\mathsf{T}} e_{\vec{p}}^k|\!|_{L^2(\Omega)}^2\leq \frac{1}{\lambda}|\!|e_{\xi}^{k-1}|\!|_{L^2(\Omega)} |\!|\vec{\alpha}^{\mathsf{T}} e_{\vec{p}}^k|\!|_{L^2(\Omega)}.
\end{align}
Then, applying the Cauchy-Schwarz inequality, we have
\begin{align*}
   \delta \frac{|\!|\vec{\alpha}^{\mathsf{T}} e_{\vec{p}}^k|\!|_{L^2(\Omega)}^2}{|\!|\vec{\alpha}|\!|^2}+ \frac{1}{\lambda} |\!|\vec{\alpha}^{\mathsf{T}} e_{\vec{p}}^k|\!|_{L^2(\Omega)}^2 &\leq \delta |\!|e_{\vec{p}}^k|\!|_{L^2(\Omega)}^2   +\frac{1}{\lambda} |\!|\vec{\alpha}^{\mathsf{T}} e_{\vec{p}}^k|\!|_{L^2(\Omega)}^2\\ \notag
   & \leq \left(  Se_{\vec{p}}^k, e_{\vec{p}}^k\right) +\frac{1}{\lambda} |\!|\vec{\alpha}^{\mathsf{T}} e_{\vec{p}}^k|\!|_{L^2(\Omega)}^2 \\ \notag
   &\leq \frac{1}{\lambda}|\!|e_{\xi}^{k-1}|\!|_{L^2(\Omega)} |\!|\vec{\alpha}^{\mathsf{T}} e_{\vec{p}}^k|\!|_{L^2(\Omega)},
\end{align*}
where $\delta=\min_{1\leq i\leq N} c_{i}$.\\
In the above inequality, dividing both sides by $|\!|\vec{\alpha}^{\mathsf{T}} e_{\vec{p}}^k|\!|_{L^2(\Omega)}$ leads to
\begin{align}
\label{decoupled_alg_ana6}
   |\!|\vec{\alpha}^{\mathsf{T}} e_{\vec{p}}^k|\!|_{L^2(\Omega)}\leq \frac{\frac{|\!|\vec{\alpha}|\!|^2}{\lambda}}{\delta+\frac{|\!|\vec{\alpha}|\!|^2}{\lambda}}|\!|e_{\xi}^{k-1}|\!|_{L^2(\Omega)}.
\end{align}
Ignoring the nonnegative term  $\frac{1}{\lambda} |\!|\vec{\alpha}^{\mathsf{T}} e_{\vec{p}}^k|\!|_{L^2(\Omega)}^2$ in  \eqref{decoupled_alg_ana8} and applying
 \eqref{decoupled_alg_ana6},
we obtain
\begin{align*}
      \delta |\!|e_{\vec{p}}^k|\!|_{L^2(\Omega)}^2\leq  \left(  Se_{\vec{p}}^k, e_{\vec{p}}^k\right)&\leq \frac{1}{\lambda}|\!|e_{\xi}^{k-1}|\!|_{L^2(\Omega)} |\!|\vec{\alpha}^{\mathsf{T}} e_{\vec{p}}^k|\!|_{L^2(\Omega)}\\ \notag
      & \leq \frac{1}{\lambda}|\!|e_{\xi}^{k-1}|\!|_{L^2(\Omega)} \frac{\frac{|\!|\vec{\alpha}|\!|^2}{\lambda}}{\delta+\frac{|\!|\vec{\alpha}|\!|^2}{\lambda}}|\!|e_{\xi}^{k-1}|\!|_{L^2(\Omega)}.
\end{align*}
By taking the square roots of both sides, we have
\begin{align} \label{decoupled_alg_ana13}
|\!|e_{\vec{p}}^k|\!|_{L^2(\Omega)} \leq \left( \frac{|\!|\vec{\alpha}|\!|^2}{\lambda \delta (\lambda \delta +|\!|\vec{\alpha}|\!|^2)} \right)^{\frac{1}{2}} |\!|e_{\xi}^{k-1}|\!|_{L^2(\Omega)}.
\end{align}
Similarly, due to $\left( \left( S+\frac{1}{\lambda} \vec{\alpha}\vec{\alpha}^{\mathsf{T}} \right)e_{\vec{p}}^k,e_{\vec{p}}^k\right) \geq 0$ and $\Delta t \left(  Be_{\vec{p}}^k, e_{\vec{p}}^k \right) \geq 0$ in  \eqref{decoupled_alg_ana3},  we derive that
\begin{align*}
\Delta t \zeta |\!|\nabla e_{\vec{p}}^k|\!|_{L^2(\Omega)}^2\leq \Delta t \left(  K\nabla e_{\vec{p}}^k,\nabla e_{\vec{p}}^k \right) &\leq \frac{1}{\lambda}\left(  \vec{\alpha}e_{\xi}^{k-1},e_{\vec{p}}^k \right) \\ \notag
& \leq \frac{1}{\lambda} |\!|e_{\xi}^{k-1}|\!| \frac{\frac{|\!|\vec{\alpha}|\!|^2}{\lambda}}{\delta+\frac{|\!|\vec{\alpha}|\!|^2}{\lambda}}|\!|e_{\xi}^{k-1}|\!|_{L^2(\Omega)},
\end{align*}
where
$\zeta=\min_{1\leq i\leq N} K_i$.\\
By taking the square roots of both sides, we obtain
\begin{align} \label{decoupled_alg_ana14}
    |\!|\nabla e_{\vec{p}}^k|\!|_{L^2(\Omega)}\leq \left(\frac{\frac{|\!|\vec{\alpha}|\!|^2}{\lambda}}{(\lambda \delta+|\!|\vec{\alpha}|\!|^2)\Delta t \zeta}\right)^{\frac{1}{2}}|\!|e_{\xi}^{k-1}|\!|_{L^2(\Omega)}.
\end{align}

Summing up \eqref{decoupled_alg_ana1} and \eqref{decoupled_alg_ana2}, there holds
\begin{align} \label{decoupled_alg_ana16}
    2\mu |\!|\epsilon(e_{\pmb{u}}^k)|\!|_{L^2(\Omega)}^2 + \frac{1}{\lambda}|\!|e_{\xi}^k|\!|_{L^2(\Omega)}^2&=\frac{1}{\lambda}(\vec{\alpha}^{\mathsf{T}} e_{\vec{p}}^k,e_{\xi}^k)
    \leq \frac{1}{\lambda}|\!|\vec{\alpha}^{\mathsf{T}} e_{\vec{p}}^k|\!|_{L^2(\Omega)} |\!|e_{\xi}^k|\!|_{L^2(\Omega)}.
\end{align}
Combining \eqref{decoupled_alg_ana16} with \eqref{decoupled_alg_ana6}, we have
\begin{align} \label{decoupled_alg_ana12}
     |\!|e_{\xi}^k|\!|_{L^2(\Omega)}\leq |\!|\vec{\alpha}^{\mathsf{T}}e_{\vec{p}}^k|\!|_{L^2(\Omega)}\leq \frac{\frac{|\!|\vec{\alpha}|\!|^2}{\lambda}}{\delta+\frac{|\!|\vec{\alpha}|\!|^2}{\lambda}}|\!|e_{\xi}^{k-1}|\!|_{L^2(\Omega)}.
\end{align}

Applying \eqref{Korn_conclusion} to \eqref{decoupled_alg_ana1} , it holds that
\begin{align*}
    2\mu |\!|\epsilon(e_{\pmb{u}}^k)|\!|_{L^2(\Omega)}^2=(e_{\xi}^k,\operatorname{div} e_{\pmb{u}}^k)&\leq |\!|e_{\xi}^k|\!|_{L^2(\Omega)}|\!|\operatorname{div} e_{\pmb{u}}^k|\!|_{L^2(\Omega)}
    \leq \sqrt{d} |\!|e_{\xi}^k|\!|_{L^2(\Omega)}|\!|\epsilon(e_{\pmb{u}}^k)|\!|_{L^2(\Omega)}.
\end{align*}
Using the first Korn inequality \cite{gu2023iterative}, it can be seen that
\begin{align} \label{decoupled_alg_ana15}
   |\!|e_{\pmb{u}}^k|\!|_{L^2(\Omega)}\leq |\!|e_{\pmb{u}}^k|\!|_{H^1(\Omega)} \leq C_K|\!|\epsilon(e_{\pmb{u}}^k)|\!|_{L^2(\Omega)} \leq
    \frac{C_K \sqrt{d}}{2 \mu}|\!|e_{\xi}^k|\!|_{L^2(\Omega)}.
\end{align}
According to \eqref{decoupled_alg_ana12} and the inequality $0<\frac{\frac{|\!|\vec{\alpha}|\!|^2}{\lambda}}{\delta+\frac{|\!|\vec{\alpha}|\!|^2}{\lambda}} <1$, we then have
  \begin{align*}
    \lim _{k\rightarrow +\infty} |\!|e_{\xi}^k|\!|_{L^2(\Omega)}=0.
\end{align*}
Following from \eqref{decoupled_alg_ana13}, \eqref{decoupled_alg_ana14} and \eqref{decoupled_alg_ana15}, we have
   \begin{align} \label{decoupled_alg_con10}
     \lim _{k\rightarrow +\infty} |\!|\nabla e_{\vec{p}}^k|\!|_{L^2(\Omega)}=
     \lim _{k\rightarrow +\infty} |\!|e_{\vec{p}}^k|\!|_{L^2(\Omega)}=
     \lim _{k\rightarrow +\infty} |\!|e_{\pmb{u}}^k|\!|_{H^1(\Omega)}=
     \lim _{k\rightarrow +\infty} |\!|e_{\pmb{u}}^k|\!|_{L^2(\Omega)}=0.
\end{align}
$\hfill\qedsymbol$

For $\min_{1\leq i \leq N} c_{i}=0$, the main conclusions are presented in the following remark.
\begin{remark} \label{Remark1}
If $\min_{1\leq i \leq N} c_{i}=0$, we will show that  \eqref{decoupled_alg_con10} also hold true. It means that
the iteratively decoupled algorithm also converges when $\min_{1\leq i \leq N} c_{i}$ degenerates to $0$.
According to \eqref{decoupled_alg_ana8}, we obtain
\begin{align} \label{decoupled_alg_ana9}
|\!|\vec{\alpha}^{\mathsf{T}} e_{\vec{p}}^k|\!|_{L^2(\Omega)}\leq |\!|e_{\xi}^{k-1}|\!|_{L^2(\Omega)}.
\end{align}
Thanks to \eqref{decoupled_alg_ana3}, \eqref{decoupled_alg_ana16}, \eqref{decoupled_alg_ana9}, and the Poincar$\rm{\acute{e}}$ inequality, we derive that
\begin{align} \label{decoupled_alg_ana17}
    |\!|e_{\xi}^k|\!|_{L^2(\Omega)}& \leq |\!|e_{\xi}^{k-1}|\!|_{L^2(\Omega)},\\
\label{decoupled_alg_ana10}
  2\mu |\!|\epsilon(e_{\pmb{u}}^k)|\!|_{L^2(\Omega)}^2 + \frac{1}{\lambda}|\!|e_{\xi}^k|\!|_{L^2(\Omega)}^2&\leq \frac{1}{\lambda}|\!|e_{\xi}^{k-1}|\!|_{L^2(\Omega)} |\!|e_{\xi}^k|\!|_{L^2(\Omega)},\\
\label{decoupled_alg_con7}
|\!|e_{\vec{p}}^k|\!|_{L^2(\Omega)} \leq C_P|\!|\nabla e_{\vec{p}}^k|\!|_{L^2(\Omega)}&\leq C_P \left( \frac{1}{\lambda \Delta t \zeta} \right)^{\frac{1}{2}} |\!|e_{\xi}^{k-1}|\!|_{L^2(\Omega)},
\end{align}
where $C_P$ is the constant in the Poincar$\rm{\acute{e}}$  inequality.
We are going to use the method of contradiction to show that the limit of $\{|\!|e_{\xi}^k|\!|_{L^2(\Omega)}\}$ is 0. If not, let us assume
$$
\lim _{k\rightarrow +\infty} |\!|e_{\xi}^k|\!|_{L^2(\Omega)}>0.
$$
From  \eqref{decoupled_alg_ana17}   and \eqref{decoupled_alg_ana10}, we know that
$$
\lim_{k\rightarrow +\infty} |\!|\epsilon(e_{\pmb{u}}^k)|\!|_{L^2(\Omega)}=0.
$$
Applying the discrete inf-sup condition, we see that
\begin{align} \label{decoupled_alg_ana11}
     \beta_0^* |\!|e_{\xi}^k|\!|_{L^2(\Omega)}\leq \sup_{\pmb{v}_h \in \pmb{V}_h} \frac{|(e_{\xi}^k,\operatorname{div} \pmb{v}_h)|}{|\!|\pmb{v}_h|\!|_{H^1(\Omega)}} =\sup_{\pmb{v}_h \in \pmb{V}_h} \frac{2\mu |(\epsilon(e_{\pmb{u}}^k),\epsilon(\pmb{v}_h))|}{|\!|\pmb{v}_h|\!|_{H^1(\Omega)}}\leq 2\mu|\!|\epsilon(e_{\pmb{u}}^k)|\!|_{L^2(\Omega)}.
\end{align}
Letting $k\rightarrow \infty$ in \eqref{decoupled_alg_ana11}, we derive  that
$$
\lim _{k\rightarrow +\infty} |\!|e_{\xi}^k|\!|_{L^2(\Omega)}\leq 0,
$$ which is a contradiction.
Therefore
\begin{align} \label{Remark2_con}
\lim _{k\rightarrow +\infty} |\!|e_{\xi}^k|\!|_{L^2(\Omega)}=0.
\end{align}
Combining \eqref{decoupled_alg_ana15} and \eqref{decoupled_alg_con7} with \eqref{Remark2_con}, it follows that \eqref{decoupled_alg_con10} holds true.
\end{remark}

\section{Numerical experiments}
In this section, we present several numerical experiments to assess the 
performance of the iteratively decoupled algorithm. 
First, accuracy tests under different physical parameter settings are performed to verify the convergence of the numerical solutions. 
Then, we conduct more physiologically realistic three-dimensional brain edema simulations.

\subsection{Accuracy tests}
Let $\Omega=[0,1]\times[0,1]$, $\Gamma_1=\{(1,y);0\leq y\leq 1\}$, $\Gamma_2=\{(x,0);0\leq x\leq 1\}$, $\Gamma_3=\{(0,y);0\leq y\leq 1\}$, and $\Gamma_4=\{(x,1);0\leq x\leq 1\}$. The terminal time is $T=0.01$. In our numerical experiments, we impose a pure Dirichlet boundary condition on this domain. We initialize the mesh with a size of $h=\frac{1}{8}$ and subsequently carry out four rounds of mesh refinement, which involves connecting the midpoints of each triangle. We use the variable $iters$ to represent the number of iterations utilized within the iteratively decoupled algorithm. In these tests, we intentionally employ relatively large time step sizes to show the efficiency of the iteratively decoupled algorithm. Specifically, we choose the values of $\Delta t$ and $iters$ such that the overall computational expense is approximately equal to that of the coupled algorithm. More precisely, we set the values of $\Delta t$ and $iters$ such that $T/ \Delta t \times iters=50$. At the end of the simulations, we compute the $H^1$ and $L^2$ norms at the final time $T$. 
All the following accuracy tests are conducted using the open-source software FreeFEM++ \cite{hecht2012new}.

The exact solutions of \eqref{three_field_form1}-\eqref{three_field_form3} are given as follows.
\begin{align*}
    \pmb{u}&=\left[\begin{array}{c}
       \text{sin}(2\pi y)(-1+\text{cos}(2\pi x))+\frac{1}{\mu+\lambda}\text{sin}(\pi x)\text{sin}(\pi y)   \\
    \text{sin}(2\pi x)(1-\text{cos}(2\pi y))+\frac{1}{\mu+\lambda}\text{sin}(\pi x) \text{sin}(\pi y)
    \end{array}
    \right]\text{sin}(t),\\
  p_1&=-\text{sin}(\pi x )\text{sin}(\pi y)\text{cos}(t),\\
    p_2&=-2 \text{sin}(\pi x)\text{sin}(\pi y)\text{cos}(t).
\end{align*}

\subsubsection{Tests for the Poisson ratio} 
In this example, we test the performance of the algorithms under different settings of the Poisson ratio $\nu$.
For other physical parameters, we set
\begin{align*}
  E=1, \alpha_1=\alpha_2=1, c_{1}=c_{2}=1, K_1=K_2=1, \beta_{12}=\beta_{21}=1.
\end{align*}

\begin{table}[ht!]
    \centering
    \caption{Convergence rate of the coupled algorithm. $\nu=0.3$ and $\Delta t=2\times 10^{-4}$.}
    \begin{tabular}{ccccc}
\hline
  $1/h$   & $L^2 \& H^1$ errors of $\pmb{u}$ & Orders & $L^2 \& H^1$ errors of $\xi$ & Orders \\
  \hline
  8 &  1.230e-03 \& 1.768e-02  &  & 3.652e-02 \& 1.083e+00  & \\
  16 & 3.013e-04 \& 4.032e-03 & 2.03 \& 2.13 & 9.105e-03 \& 5.506e-01 & 2.00 \& 0.976\\
  32& 7.536e-05 \& 9.421e-04 & 2.00 \& 2.10 & 2.269e-03 \& 2.760e-01 & 2.00 \& 0.996\\
  64& 1.890e-05 \& 2.257e-04 & 2.00 \& 2.06 & 5.670e-04 \& 1.381e-01 & 2.00 \& 0.999\\
  128& 4.766e-06 \& 5.523e-05 & 1.99 \& 2.03 & 1.423e-04 \& 6.908e-02 & 1.99 \& 1.000\\
  \hline

  $1/h$ & $L^2 \& H^1$ errors of $p_1$ & Orders & $L^2 \& H^1$ errors of $p_2$&Orders\\
  \hline
  8 &1.432e-02 \& 3.581e-01 &   & 2.851e-02 \&  7.161e-01 & \\
  16& 3.681e-03 \& 1.816e-01 &1.96 \& 0.98 &7.342e-03 \& 3.633e-01 & 1.96 \& 0.98\\
  32 &  9.354e-04 \& 9.134e-02 &1.98 \& 0.99 & 1.868e-03 \& 1.827e-01 &1.97 \& 0.99\\
  64& 2.403e-04 \& 4.576e-02 &1.96 \& 1.00 &4.809e-04 \& 9.153e-02  &1.96 \& 1.00\\
  128& 6.586e-05 \& 2.290e-02 & 1.87 \&  1.00 & 1.327e-04 \& 4.579e-02& 1.86 \& 1.00\\
  \hline
\end{tabular}
    \label{table 1}
\end{table}

In Tables \ref{table 1} through \ref{table 3}, we present the numerical results for both the coupled algorithm and the iteratively decoupled algorithm in the case of Poisson's ratio $\nu=0.3$. Table \ref{table 2} provides the results derived from the iteratively decoupled algorithm with 10 iterations, while Table \ref{table 3} showcases the results for 20 iterations. The numerical results presented in Tables \ref{table 1} to \ref{table 3} demonstrate that the energy-norm errors of all variables converge to their optimal orders. Furthermore, the results presented in Table \ref{table 3} illustrate that increasing the number of iterations enhances the accuracy of the iteratively decoupled algorithm.

\begin{table}[ht!]
    \centering
    \caption{Convergence rate of the iteratively decoupled algorithm. $\nu=0.3$, $\Delta t=2\times 10^{-3}$, and $iters=10$.}
   \begin{tabular}{ccccc}
\hline
  $1/h$   & $L^2 \& H^1$ errors of $\pmb{u}$ &Orders& $L^2 \& H^1$ errors of $\xi$ & Orders\\
  \hline
  8 & 1.229e-03 \& 1.765e-02 &  & 3.667e-02 \& 1.084e+00  & \\
  16 & 3.011e-04 \& 4.012e-03 & 2.03 \& 2.14 & 9.146e-03 \& 5.507e-01 & 2.00 \& 0.977 \\
  32& 7.577e-05 \& 9.367e-04 & 1.99 \& 2.10 & 2.283e-03 \& 2.760e-01 & 2.00 \& 0.996\\
  64& 1.961e-05 \& 2.260e-04 & 1.95 \& 2.05 & 5.738e-04 \& 1.381e-01 & 1.99 \& 0.999\\
  128& 5.921e-06 \& 5.813e-05 & 1.73 \& 1.96 & 1.478e-04 \& 6.908e-02 & 1.96 \& 1.000\\
  \hline

  $1/h$ & $L^2 \& H^1$ errors of $p_1$ &Orders & $L^2 \& H^1$ errors of $p_2$&Orders\\
  \hline
  8 & 1.200e-02 \& 3.609e-01 &   & 2.625e-02 \& 7.203e-01 & \\
  16& 3.036e-03 \& 1.827e-01 &1.98 \& 0.98 & 6.662e-03 \& 3.652e-01 & 1.98 \& 0.98\\
  32 &  7.626e-04 \& 9.182e-02 & 1.99 \& 0.99 & 1.678e-03 \& 1.836e-01 & 1.99 \& 0.99\\
  64& 1.920e-04 \& 4.599e-02 & 1.99 \& 1.00 & 4.239e-04 \& 9.198e-02 & 1.98 \& 1.00\\
  128& 4.950e-05 \& 2.301e-02 & 1.96 \& 1.00 & 1.100e-04 \& 4.602e-02 & 1.95 \& 1.00\\
  \hline
\end{tabular}
    \label{table 2}
\end{table}

\begin{table}[ht!]
    \centering
    \caption{Convergence rate of the iteratively decoupled algorithm. $\nu=0.3$, $\Delta t=4\times 10^{-3}$, and $iters=20$.}
   \begin{tabular}{ccccc}
\hline
  $1/h$   & $L^2 \& H^1$ errors of $\pmb{u}$ &Orders& $L^2 \& H^1$ errors of $\xi$&Orders\\
  \hline
  8 & 1.184e-03 \& 1.734e-02 &  & 3.656e-02 \& 1.084e+00  & \\
  16 & 2.879e-04 \& 3.877e-03 & 2.04 \& 2.16 & 9.114e-03 \& 5.508e-01 & 2.00 \& 0.977\\
  32& 7.230e-05 \& 8.955e-04 & 1.99 \& 2.11 & 2.279e-03 \& 2.760e-01 & 2.00 \& 0.997 \\
  64& 1.858e-05 \& 2.150e-04 & 1.96 \& 2.06 & 5.782e-04 \& 1.381e-01 & 1.98 \& 0.999 \\
  128& 5.148e-06 \& 5.469e-05 & 1.85 \& 1.97 & 1.542e-04 \& 6.908e-02 & 1.91 \& 1.000 \\
  \hline
  $1/h$ & $L^2 \& H^1$ errors of $p_1$ &Orders & $L^2 \& H^1$ errors of $p_2$ &Orders\\
  \hline
  8 & 1.192e-02 \& 3.613e-01 &   & 2.613e-02 \& 7.206e-01 & \\
  16& 3.004e-03 \& 1.828e-01 &1.99 \& 0.98 & 6.616e-03 \& 3.652e-01 &1.98 \& 0.98\\
  32 &  7.554e-04 \& 9.183e-02 & 1.99 \& 0.99  & 1.668e-03 \& 1.836e-01 & 1.99 \& 0.99\\
  64& 1.919e-04 \& 4.600e-02 & 1.98 \& 1.00 & 4.248e-04 \& 9.199e-02& 1.97 \& 1.00\\
  128& 5.088e-05 \& 2.301e-02 & 1.91 \& 1.00 & 1.136e-04 \& 4.602e-02 & 1.90 \& 1.00\\
  \hline
\end{tabular}
    \label{table 3}
\end{table}

Tables \ref{table 1} to \ref{table 3} pertain to the scenario in which the poroelastic material is compressible. In contrast, Tables \ref{table 4} to \ref{table 6} maintain constant physical parameters while specifying a Poisson ratio of $\nu=0.49999.$ This specific Poisson ratio reflects the nearly incompressible behavior of the poroelastic material. Both Table \ref{table 4} and Tables \ref{table 5} to \ref{table 6} present results derived from the coupled algorithm and the iteratively decoupled algorithm, respectively, each with different iteration counts. Notably, as the Poisson ratio $\nu$ approaches 0.5, indicating that the elastic material is almost incompressible, both algorithms exhibit enhanced numerical error rates in comparison to the case when $\nu=0.3$.
Upon examining Tables \ref{table 4} to \ref{table 6}, it becomes evident that, in both algorithms, the $L^2$ norm errors and energy norm errors for each variable demonstrate optimal convergence rates. A comparative analysis of the numerical outcomes presented in Tables \ref{table 4} and \ref{table 6} underscores the superior error accuracy achieved by our proposed iteratively decoupled algorithm.

\begin{table}[ht!]
    \centering
    \caption{Convergence rate of the coupled algorithm. $\nu=0.49999$ and $\Delta t=2\times 10^{-4}$.}
    \begin{tabular}{ccccc}
\hline
  $1/h$   & $L^2 \& H^1$ errors of $\pmb{u}$ & Orders & $L^2 \& H^1$ errors of $\xi$ & Orders \\
  \hline
  8 &  4.109e-04 \& 1.668e-02  &  & 3.945e-02 \& 1.090e+00  & \\
  16 & 3.574e-05 \& 3.081e-03 & 3.52 \& 2.44 & 9.785e-03 \& 5.607e-01 & 2.01 \& 0.959\\
  32& 3.262e-06 \& 5.696e-04 & 3.45 \& 2.44 & 2.431e-03 \& 2.779e-01 & 2.01 \& 1.013\\
  64& 2.975e-07 \& 1.068e-04 & 3.45 \& 2.42 & 6.065e-04 \& 1.386e-01 & 2.00 \& 1.004\\
  128& 2.776e-08 \& 2.076e-05 & 3.42 \& 2.36 & 1.521e-04 \& 6.919e-02 & 2.00 \& 1.002\\
  \hline
  $1/h$ & $L^2 \& H^1$ errors of $p_1$ & Orders & $L^2 \& H^1$ errors of $p_2$&Orders\\
  \hline
  8 & 1.469e-02 \& 3.599e-01 &   & 2.911e-02 \&  7.197e-01 & \\
  16& 3.743e-03 \& 1.825e-01 &1.97 \& 0.98 &7.414e-03 \& 3.651e-01 & 1.97 \& 0.98\\
  32 &  9.426e-04 \& 9.180e-02 &1.99 \& 0.99 & 1.867e-03 \& 1.836e-01 &1.99 \& 0.99\\
  64& 2.365e-04 \& 4.599e-02 &1.99 \& 1.00 & 4.684e-04 \& 9.198e-02  &1.99 \& 1.00\\
  128& 5.945e-05 \& 2.301e-02 & 1.99 \& 1.00 & 1.177e-04 \& 4.602e-02& 1.99 \& 1.00\\
  \hline
\end{tabular}
    \label{table 4}
\end{table}

\begin{table}[ht!]
    \centering
    \caption{Convergence rate of the iteratively decoupled algorithm. $\nu=0.49999$, $\Delta t=2\times 10^{-3}$, and $iters=10$.}
   \begin{tabular}{ccccc}
\hline
  $1/h$   & $L^2 \& H^1$ errors of $\pmb{u}$ &Orders& $L^2 \& H^1$ errors of $\xi$ & Orders\\
  \hline
  8 & 4.109e-04 \& 1.668e-02 &  & 3.942e-02 \& 1.090e+00  & \\
  16 & 3.574e-05 \& 3.081e-03 & 3.52 \& 2.44 & 9.783e-03 \& 5.607e-01 & 2.01 \& 0.959 \\
  32& 3.262e-06 \& 5.696e-04 & 3.45 \& 2.44 & 2.437e-03 \& 2.779e-01 & 2.01 \& 1.013\\
  64& 2.975e-07 \& 1.068e-04 & 3.45 \& 2.42 & 6.144e-04 \& 1.386e-01 & 1.99 \& 1.004\\
  128& 2.776e-08 \& 2.076e-05 & 3.42 \& 2.36 & 1.607e-04 \& 6.919e-02 & 1.94 \& 1.002\\
  \hline

  $1/h$ & $L^2 \& H^1$ errors of $p_1$ &Orders & $L^2 \& H^1$ errors of $p_2$&Orders\\
  \hline
  8 & 1.469e-02 \& 3.599e-01 &   & 2.909e-02 \& 7.198e-01 & \\
  16& 3.743e-03 \& 1.825e-01 &1.97 \& 0.98 & 7.413e-03 \& 3.651e-01 & 1.97 \& 0.98\\
  32 &  9.449e-04 \& 9.180e-02 & 1.99 \& 0.99 & 1.872e-03 \& 1.836e-01 & 1.99 \& 0.99\\
  64& 2.396e-04 \& 4.599e-02 & 1.98 \& 1.00 & 4.746e-04 \& 9.198e-02 & 1.98 \& 1.00\\
  128& 6.281e-05 \& 2.301e-02 & 1.93 \& 1.00 & 1.244e-04 \& 4.602e-02 & 1.93 \& 1.00\\
  \hline
\end{tabular}
    \label{table 5}
\end{table}

\newpage

\begin{table}[ht!]
    \centering
    \caption{Convergence rate of the iteratively decoupled algorithm. $\nu=0.49999$, $\Delta t=4\times 10^{-3}$, and $iters=20$.}
   \begin{tabular}{ccccc}
\hline
  $1/h$   & $L^2 \& H^1$ errors of $\pmb{u}$ &Orders& $L^2 \& H^1$ errors of $\xi$ & Orders\\
  \hline
  8 & 4.096e-04 \& 1.657e-02 &  & 3.927e-02 \& 1.090e+00  & \\
  16 & 3.552e-05 \& 3.046e-03 & 3.53 \& 2.44 & 9.738e-03 \& 5.606e-01 & 2.01 \& 0.959 \\
  32& 3.224e-06 \& 5.575e-04 & 3.46 \& 2.45 & 2.429e-03 \& 2.779e-01 & 2.00 \& 1.012\\
  64& 2.911e-07 \& 1.027e-04 & 3.47 \& 2.44 & 6.168e-04 \& 1.386e-01 & 1.98 \& 1.004\\
  128& 2.669e-08 \& 1.944e-05 & 3.45 \& 2.40 & 1.657e-04 \& 6.918e-02 & 1.90 \& 1.002\\
  \hline

  $1/h$ & $L^2 \& H^1$ errors of $p_1$ &Orders & $L^2 \& H^1$ errors of $p_2$&Orders\\
  \hline
  8 & 1.462e-02 \& 3.600e-01 &   & 2.903e-02 \& 7.199e-01 & \\
  16& 3.723e-03 \& 1.826e-01 & 1.97 \& 0.98 & 7.388e-03 \& 3.651e-01 & 1.97 \& 0.98\\
  32 &  9.410e-04 \& 9.180e-02 & 1.98 \& 0.99 & 1.867e-03 \& 1.836e-01 & 1.98 \& 0.99\\
  64& 2.403e-04 \& 4.599e-02 & 1.97 \& 1.00 & 4.768e-04 \& 9.198e-02 & 1.97 \& 1.00\\
  128& 6.476e-05 \& 2.301e-02 & 1.89 \& 1.00 & 1.285e-04 \& 4.602e-02 & 1.89 \& 1.00\\
  \hline
\end{tabular}
    \label{table 6}
\end{table}


\subsubsection{Tests for the hydraulic conductivity} 
In this example, we conduct experiments to assess the performance of the iteratively decoupled algorithm by varying the values of the parameter $K_i$. It is worth noting that we previously examined the case where $K_1=K_2=1$ in our earlier tests, but now we set $K_1=K_2=10^{-6}$. For the other physical parameters, we maintain the following values:
\begin{align*}
E=1, \nu=0.3, \alpha_1=\alpha_2=1, c_{1}=c_{2}=1, \beta_{12}=\beta_{21}=1.
\end{align*}

The results presented in Table \ref{table 8} through Table \ref{table 9} illustrate that increasing the number of iterations results in improved accuracy when using the iteratively decoupled algorithm. Upon analyzing the numerical results in Table \ref{table 3} and Table \ref{table 9}, it becomes apparent that the energy-norm errors for all variables converge to their optimal orders. This suggests that the accuracy of the iteratively decoupled algorithm remains robust and is minimally affected by variations in hydraulic conductivity.

\begin{table}[ht!]
    \centering
    \caption{Convergence rate of the coupled algorithm. $K_1=K_2=10^{-6}$ and $\Delta t=2\times 10^{-4}$.}
    \begin{tabular}{ccccc}
\hline
  $1/h$   & $L^2 \& H^1$ errors of $\pmb{u}$ & Orders & $L^2 \& H^1$ errors of $\xi$ & Orders \\
  \hline
  8 &  1.021e-03 \& 1.674e-02  &  & 3.674e-02 \& 1.092e+00  & \\
  16 & 2.283e-04 \& 3.321e-03 & 2.16 \& 2.33 & 9.210e-03 \& 5.561e-01 & 2.00 \& 0.973\\
  32& 5.444e-05 \& 6.833e-04 & 2.07 \& 2.28 & 2.326e-03 \& 2.772e-01 & 1.99 \& 1.004\\
  64& 1.334e-05 \& 1.477e-04 & 2.03 \& 2.21 & 5.852e-04 \& 1.384e-01 & 1.99 \& 1.002\\
  128& 3.325e-06 \& 3.369e-05 & 2.00 \& 2.13 & 1.475e-04 \& 6.914e-02 & 1.99 \& 1.001\\
  \hline
  $1/h$ & $L^2 \& H^1$ errors of $p_1$ & Orders & $L^2 \& H^1$ errors of $p_2$&Orders\\
  \hline
  8 &1.210e-02 \& 3.662e-01 &   & 2.584e-02 \&  7.288e-01 & \\
  16& 3.037e-03 \& 1.863e-01 &1.99 \& 0.98 &6.453e-03 \& 3.690e-01 & 2.00 \& 0.98\\
  32 &  7.711e-04 \& 9.268e-02 &1.98 \& 1.01 & 1.622e-03 \& 1.844e-01 &1.99 \& 1.00\\
  64& 1.945e-04 \& 4.620e-02 &1.99 \& 1.00 &4.069e-04 \& 9.214e-02  &1.99 \& 1.00\\
  128& 4.911e-05 \& 2.306e-02 & 1.99 \& 1.00 & 1.024e-04 \& 4.605e-02& 1.99 \& 1.00\\
  \hline
\end{tabular}
    \label{table 7}
\end{table}

\begin{table}[ht!]
    \centering
    \caption{Convergence rate of the iteratively decoupled algorithm. $K_1=K_2=10^{-6}$, $\Delta t=2\times 10^{-3}$, and $iters=10$.}
   \begin{tabular}{ccccc}
\hline
  $1/h$   & $L^2 \& H^1$ errors of $\pmb{u}$ &Orders& $L^2 \& H^1$ errors of $\xi$ & Orders\\
  \hline
  8 & 1.025e-03 \& 1.676e-02 &  & 3.688e-02 \& 1.092e+00  & \\
  16 & 2.301e-04 \& 3.326e-03 & 2.16 \& 2.33 & 9.248e-03 \& 5.561e-01 & 2.00 \& 0.973 \\
  32& 5.618e-05 \& 6.874e-04 & 2.03 \& 2.27 & 2.337e-03 \& 2.772e-01 & 1.98 \& 1.004\\
  64& 1.588e-05 \& 1.543e-04 & 1.82 \& 2.16 & 5.926e-04 \& 1.384e-01 & 1.98 \& 1.002\\
  128& 7.569e-06 \& 4.856e-05 & 1.07 \& 1.67 & 1.557e-04 \& 6.915e-02 & 1.93 \& 1.001\\
  \hline

  $1/h$ & $L^2 \& H^1$ errors of $p_1$ &Orders & $L^2 \& H^1$ errors of $p_2$&Orders\\
  \hline
  8 & 1.217e-02 \& 3.663e-01 &   & 2.592e-02 \& 7.288e-01 & \\
  16& 3.054e-03 \& 1.863e-01 & 1.99 \& 0.98 & 6.474e-03 \& 3.690e-01 & 2.00 \& 0.98\\
  32 &  7.753e-04 \& 9.269e-02 & 1.98 \& 1.01 & 1.630e-03 \& 1.844e-01 & 1.99 \& 1.00\\
  64& 1.971e-04 \& 4.621e-02 & 1.98 \& 1.00 & 4.130e-04 \& 9.215e-02 & 1.98 \& 1.00\\
  128& 5.370e-05 \& 2.308e-02 & 1.88 \& 1.00 & 1.093e-04 \& 4.606e-02 & 1.92 \& 1.00\\
  \hline
\end{tabular}
    \label{table 8}
\end{table}

\begin{table}[ht!]
    \centering
    \caption{Convergence rate of the iteratively decoupled algorithm. $K_1=K_2=10^{-6}$, $\Delta t=4\times 10^{-3}$, and $iters=20$.}
   \begin{tabular}{ccccc}
\hline
  $1/h$   & $L^2 \& H^1$ errors of $\pmb{u}$ &Orders& $L^2 \& H^1$ errors of $\xi$&Orders\\
  \hline
  8 & 1.021e-03 \& 1.664e-02 &  & 3.675e-02 \& 1.091e+00  & \\
  16 & 2.289e-04 \& 3.291e-03 & 2.16 \& 2.34 & 9.226e-03 \& 5.560e-01 & 1.99 \& 0.973\\
  32& 5.500e-05 \& 6.757e-04 & 2.06 \& 2.28 & 2.341e-03 \& 2.772e-01 & 1.98 \& 1.004 \\
  64& 1.389e-05 \& 1.474e-04 & 1.99 \& 2.20 & 6.010e-04 \& 1.384e-01 & 1.96 \& 1.002 \\
  128& 3.871e-06 \& 3.579e-05 & 1.84 \& 2.04 & 1.636e-04 \& 6.914e-02 & 1.88 \& 1.001 \\
  \hline
  $1/h$ & $L^2 \& H^1$ errors of $p_1$ &Orders & $L^2 \& H^1$ errors of $p_2$ &Orders\\
  \hline
  8 & 1.209e-02 \& 3.663e-01 &   & 2.588e-02 \& 7.288e-01 & \\
  16& 3.038e-03 \& 1.863e-01 &1.99 \& 0.98 & 6.470e-03 \& 3.690e-01 &2.00 \& 0.98\\
  32 &  7.750e-04 \& 9.269e-02 & 1.97 \& 1.01  & 1.634e-03 \& 1.844e-01 & 1.99 \& 1.00\\
  64& 1.994e-04 \& 4.621e-02 & 1.96 \& 1.00 & 4.184e-04 \& 9.215e-02& 1.97 \& 1.00\\
  128& 5.432e-05 \& 2.307e-02 & 1.88 \& 1.00 & 1.139e-04 \& 4.606e-02 & 1.88 \& 1.00\\
  \hline
\end{tabular}
    \label{table 9}
\end{table}

\newpage

\subsubsection{Tests for the storage coefficients} 
In this example, we evaluate the influence of the specific storage coefficient $c_i$ on the iteratively decoupled algorithm. We now let $c_1=c_2=0$. For other physical parameters, we set
\begin{align*}
  E=1, \nu=0.3, \alpha_1=\alpha_2=1,  K_1=K_2=1, \beta_{12}=\beta_{21}=1.
\end{align*}

Comparing the results in Table \ref{table 2} with those in Table \ref{table 11}, it is evident that in cases where $c_1=c_2=0$, the error rates for all variables exhibit a slight degradation when utilizing the iteratively decoupled algorithm. However, a closer examination of Table \ref{table 12} reveals that as the number of iterations increases to 20, the error rates for all variables in the iteratively decoupled algorithm tend to decrease, and the energy-norm errors attain their optimal convergence rates. These experimental findings corroborate the assertion made in Remark \ref{Remark1}.  

\begin{table}[ht!]
    \centering
    \caption{Convergence rate of the coupled algorithm. $c_1=c_2=0$ and $\Delta t=2\times 10^{-4}$.}
    \begin{tabular}{ccccc}
\hline
  $1/h$   & $L^2 \& H^1$ errors of $\pmb{u}$ & Orders & $L^2 \& H^1$ errors of $\xi$ & Orders \\
  \hline
  8 &  8.027e-04 \& 1.657e-02  &  & 3.108e-02 \& 1.084e+00  & \\
  16 & 1.860e-04 \& 3.729e-03 & 2.11 \& 2.15 & 7.778e-03 \& 5.506e-01 & 2.00 \& 0.977\\
  32& 4.606e-05 \& 8.620e-04 & 2.01 \& 2.11 & 1.941e-03 \& 2.760e-01 & 2.00 \& 0.996\\
  64& 1.152e-05 \& 2.047e-04 & 2.00 \& 2.07 & 4.845e-04 \& 1.381e-01 & 2.00 \& 0.999\\
  128& 2.883e-06 \& 4.967e-05 & 2.00 \& 2.04 & 1.210e-04 \& 6.908e-02 & 2.00 \& 1.000\\
  \hline
  $1/h$ & $L^2 \& H^1$ errors of $p_1$ & Orders & $L^2 \& H^1$ errors of $p_2$&Orders\\
  \hline
  8 &9.448e-03 \& 3.622e-01 &   & 2.296e-02 \&  7.210e-01 & \\
  16& 2.417e-03 \& 1.829e-01 &1.97 \& 0.99 &5.890e-03 \& 3.652e-01 & 1.96 \& 0.98\\
  32 &  6.097e-04 \& 9.184e-02 &1.99 \& 0.99 & 1.488e-03 \& 1.836e-01 &1.98 \& 0.99\\
  64& 1.529e-04 \& 4.600e-02 &2.00 \& 1.00 &3.734e-04 \& 9.198e-02  &1.99 \& 1.00\\
  128& 3.827e-05 \& 2.301e-02 & 2.00 \& 1.00 & 9.345e-05 \& 4.602e-02& 2.00 \& 1.00\\
  \hline
\end{tabular}
    \label{table 10}
\end{table}

\begin{table}[ht!]
    \centering
    \caption{Convergence rate of the iteratively decoupled algorithm. $c_1=c_2=0$, $\Delta t=2\times 10^{-3}$, and $iters=10$.}
   \begin{tabular}{ccccc}
\hline
  $1/h$   & $L^2 \& H^1$ errors of $\pmb{u}$ &Orders& $L^2 \& H^1$ errors of $\xi$ & Orders\\
  \hline
  8 & 8.829e-04 \& 1.673e-02 &  & 3.264e-02 \& 1.084e+00  & \\
  16 & 2.080e-04 \& 3.760e-03 & 2.09 \& 2.15 & 8.127e-03 \& 5.507e-01 & 2.01 \& 0.977 \\
  32& 5.719e-05 \& 8.770e-04 & 1.86 \& 2.10 & 2.011e-03 \& 2.760e-01 & 2.01 \& 0.996\\
  64& 2.944e-05 \& 2.442e-04 & 0.96 \& 1.84 & 4.972e-04 \& 1.381e-01 & 2.02 \& 0.999\\
  128& 2.702e-05 \& 1.421e-04 & 0.12 \& 0.78 & 1.529e-04 \& 6.909e-02 & 1.70 \& 1.000\\
  \hline
  $1/h$ & $L^2 \& H^1$ errors of $p_1$ &Orders & $L^2 \& H^1$ errors of $p_2$&Orders\\
  \hline
  8 & 1.011e-02 \& 3.618e-01 &   & 2.383e-02 \& 7.208e-01 & \\
  16& 2.554e-03 \& 1.828e-01 &1.99 \& 0.98 & 6.081e-03 \& 3.652e-01 & 1.97 \& 0.98\\
  32 &  6.388e-04 \& 9.184e-02 & 2.00 \& 0.99 & 1.527e-03 \& 1.836e-01 & 1.99 \& 0.99\\
  64& 1.722e-04 \& 4.600e-02 & 1.89 \& 1.00 & 3.834e-04 \& 9.199e-02 & 1.99 \& 1.00\\
  128& 8.820e-05 \& 2.302e-02 & 0.97 \& 1.00 & 1.191e-04 \& 4.603e-02 & 1.69 \& 1.00\\
  \hline
\end{tabular}
    \label{table 11}
\end{table}

\begin{table}[ht!]
    \centering
    \caption{Convergence rate of the iteratively decoupled algorithm. $c_1=c_2=0$, $\Delta t=4\times 10^{-3}$, and $iters=20$.}
   \begin{tabular}{ccccc}
\hline
  $1/h$   & $L^2 \& H^1$ errors of $\pmb{u}$ &Orders& $L^2 \& H^1$ errors of $\xi$&Orders\\
  \hline
  8 & 7.340e-04 \& 1.620e-02 &  & 3.063e-02 \& 1.084e+00  & \\
  16 & 1.643e-04 \& 3.555e-03 & 2.16 \& 2.19 & 7.616e-03 \& 5.508e-01 & 2.01 \& 0.977\\
  32& 4.024e-05 \& 8.060e-04 & 2.03 \& 2.14 & 1.890e-03 \& 2.760e-01 & 2.01 \& 0.997 \\
  64& 9.963e-06 \& 1.890e-04 & 2.01 \& 2.09 & 4.675e-04 \& 1.381e-01 & 2.02 \& 0.999 \\
  128& 2.411e-06 \& 4.533e-05 & 2.05 \& 2.06 & 1.130e-04 \& 6.908e-02 & 2.05 \& 1.000 \\
  \hline

  $1/h$ & $L^2 \& H^1$ errors of $p_1$ &Orders & $L^2 \& H^1$ errors of $p_2$ &Orders\\
  \hline
  8 & 9.229e-03 \& 3.627e-01 &   & 2.263e-02 \& 7.213e-01 & \\
  16& 2.323e-03 \& 1.830e-01 &1.99 \& 0.99 & 5.769e-03 \& 3.653e-01 &1.97 \& 0.98\\
  32 &  5.801e-04 \& 9.186e-02 & 2.00 \& 0.99  & 1.451e-03 \& 1.836e-01 & 1.99 \& 0.99\\
  64& 1.437e-04 \& 4.600e-02 & 2.01 \& 1.00 & 3.619e-04 \& 9.199e-02& 2.00 \& 1.00\\
  128& 3.457e-05 \& 2.301e-02 & 2.06 \& 1.00 & 8.866e-05 \& 4.602e-02 & 2.03 \& 1.00\\
  \hline
\end{tabular}
    \label{table 12}
\end{table}

\subsection{Brain edema simulations}
In this work, we employ a four-network poroelastic model \cite{lee2019mixed,kraus2023hybridized} to simulate fluid dynamics and pressure distributions in brain tissue. The model characterizes four distinct vascular compartments through their respective pressures: 
  $p_1$: Interstitial fluid pressure in extracellular spaces;
    $p_2$: Arterial blood pressure; 
 $p_3$: Venous blood pressure;
 $p_4$: Capillary blood pressure.
The corresponding material parameters for these networks are specified in Table \ref{tab:brain_parameters}.

The brain geometry $\Omega$ is defined by Version 2 of the Colin 27 Adult Brain Atlas 
	FEM mesh \cite{fang2010mesh} and we use the same mesh as in \cite{lee2019mixed}, which consists of 99605 elements and 29037 vertices. The surface $\partial\Omega$ is partitioned into 
	two distinct boundaries: the outer boundary $\Gamma_{s}$ corresponding to the skull 
	and the inner boundary $\Gamma_{v}$ corresponding to the ventricles. 
    Considering the influence of the heartbeat, the boundary conditions 
	are given as follows (all values are given in mmHg, 1 mmHg $\approx$ 133.32 Pa):

	\begin{alignat*}{5}
		&p_{1}& &= 5.0 + 2.0\sin{(2\pi t)}& \quad &\textup{on}\,\, \Gamma_{s}, & 
		\qquad &p_{1} = 5.0 + 2.012\sin{(2\pi t)}& \quad &\textup{on}\,\, \Gamma_{v}, \\
		&p_{2}& &= 70.0 + 10.0\sin{(2\pi t)}&  &\textup{on}\,\, \Gamma_{s},&
		&\nabla p_{2} \cdot \bm{n} = 0&  &\textup{on}\,\, \Gamma_{v}, \\
		&p_{3}& &= 6.0& &\textup{on}\,\, \Gamma_{s}, &
		&p_{3} = 6.0& &\textup{on}\,\, \Gamma_{v}, \\
		&\nabla & &p_{4} \cdot \bm{n} = 0& &\textup{on}\,\, \Gamma_{s}, &
		&\nabla p_{4} \cdot \bm{n} = 0 & &\textup{on}\,\, \Gamma_{v}, \\
		&\bm{u}& &= \bm{0}& &\textup{on}\,\, \Gamma_{s}, &
		&(\sigma(\bm{u}) - (\vec{\alpha}^{\textup{T}} \vec{p})\bm{I}) \cdot \bm{n} = s\,\bm{n}&
		&\textup{on}\,\, \Gamma_{v}.
	\end{alignat*}
	Here, $s$ is defined by
	\begin{equation*}
		s = - \sum_{j = 1}^{4} \alpha_{j}p_{j, V},
	\end{equation*}
    where $p_{j, V}$ are the prescribed pressures on the ventricles given by:
	\begin{equation*}
		p_{1, V} = 5.0 + 2.012\sin{(2\pi t)}, \quad p_{2, V} = 70.0 + 10.0\sin{(2\pi t)},
	\end{equation*}
	\begin{equation*}
		p_{3, V} = 6.0, \quad p_{4, V} = 38.0.
	\end{equation*}
    The initial conditions are given by
    \begin{equation*}
        \bm{u}_{0} = \bm{0}, \quad p_{1, 0} = 5.0 \quad p_{2, 0} = 70.0, \quad
        p_{3, 0} = 6.0 \quad p_{4, 0} = 38.0.
    \end{equation*}

	\begin{figure}[!htbp]
		\centering
		\includegraphics[scale=0.25]{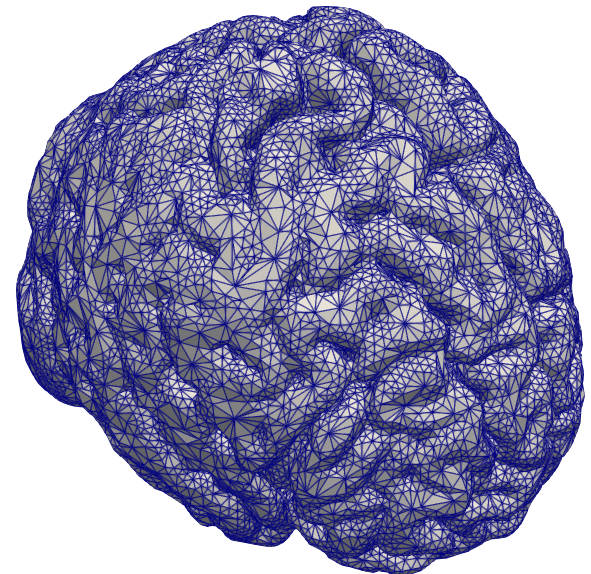}
		\caption{A 3D brain mesh.}
		\label{fig:brain_mesh}
	\end{figure}

	\begin{table}[!htbp]
		\centering
        \caption{Reference values of parameters.}
		\begin{tabular}{ccc}
			\hline
			\textbf{Parameter} & \textbf{Value} & \textbf{Unit}     \\
			\hline
			$\nu$               & 0.4999             &              \\
			$E$                 & 1500               & Pa           \\
			$c_{1}$             & $3.9\times 10^{-4}$& Pa$^{-1}$    \\
			$c_{2}, c_{4}$      & $2.9\times 10^{-4}$& Pa$^{-1}$    \\
			$c_{3}$             & $1.5\times 10^{-5}$& Pa$^{-1}$    \\
			$\alpha_{1}$        & 0.49               &              \\
			$\alpha_{2}, \alpha_{4}$ & 0.25          &              \\
			$\alpha_{3}$        & 0.01               &              \\
			$K_{1}$             & $1.57\times 10^{-5}$ & mm$^{2}$Pa$^{-1}$s$^{-1}$ \\
			$K_{2}, K_{3}, K_{4}$ & $3.75\times 10^{-6}$ & mm$^{2}$Pa$^{-1}$s$^{-1}$ \\
			$\beta_{13}, \beta_{14}, \beta_{24}, \beta_{34}$ & $1.0\times 10^{-6}$ & Pa$^{-1}$s$^{-1}$ \\
			$\beta_{12}, \beta_{23}$ & 0.0          & Pa$^{-1}$s$^{-1}$ \\
			\hline
		\end{tabular}
		\label{tab:brain_parameters}
	\end{table}
We implement our brain edema simulation algorithms using the open-source finite element package FEniCS~\cite{logg2012automated}. To evaluate computational efficiency, we compare the CPU time required by both the iteratively decoupled and coupled algorithms to reach the final simulation time $T = 3.0\,\text{s}$. 

The numerical experiments employ the following time step configurations:
\begin{itemize}
    \item Iteratively decoupled algorithm: $\Delta t_1 = 0.0625\,\text{s}$ with $5$ iterations per step
    \item Coupled algorithm: $\Delta t_2 = 0.0125\,\text{s}$
\end{itemize}
On the initial coarse mesh, the decoupled algorithm completes in $778\,\text{s}$ compared to $1387\,\text{s}$ for the coupled approach. We further validate these results on a refined mesh containing 249,361 elements and 55,066 vertices, where the decoupled algorithm requires $3189\,\text{s}$ versus $5874\,\text{s}$ for the coupled version. These results consistently demonstrate the superior computational efficiency of our iteratively decoupled approach.

\begin{figure}[!htbp]
		\centering
		\begin{subfigure}[b]{.45\linewidth}
			\centering
			\includegraphics[width=6.0cm,height=5.0cm]{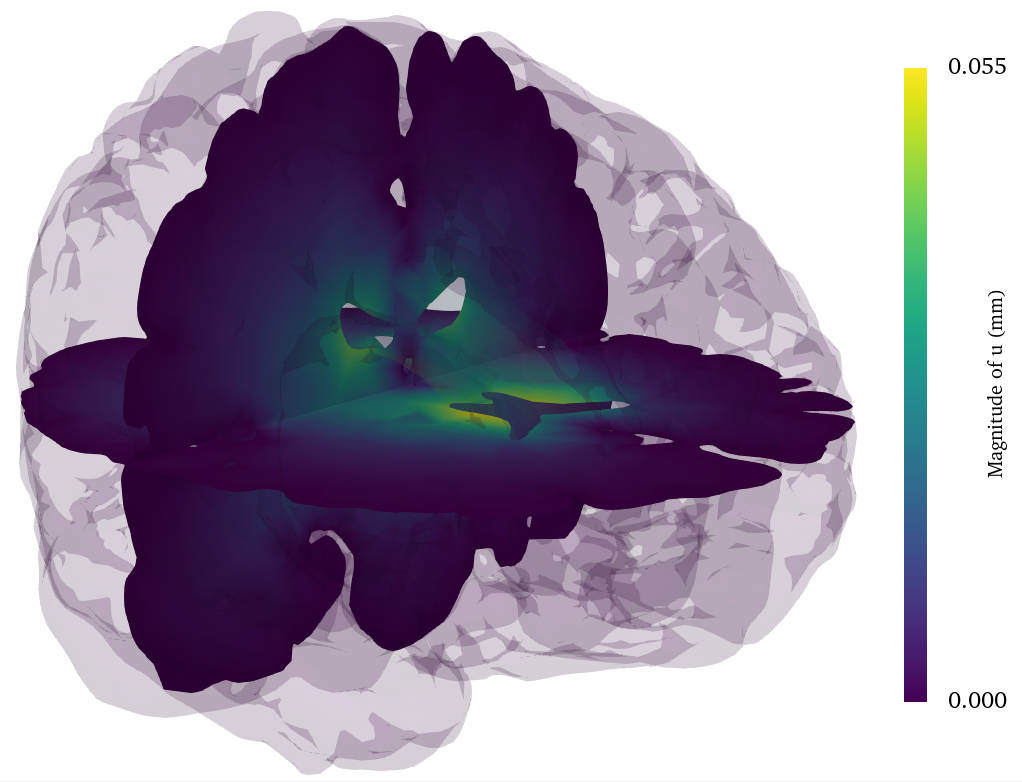}
			\caption{the iteratively decoupled algorithm}
		\end{subfigure}
        \qquad
		\begin{subfigure}[b]{.45\linewidth}
			\centering
			\includegraphics[width=6.0cm,height=5.0cm]{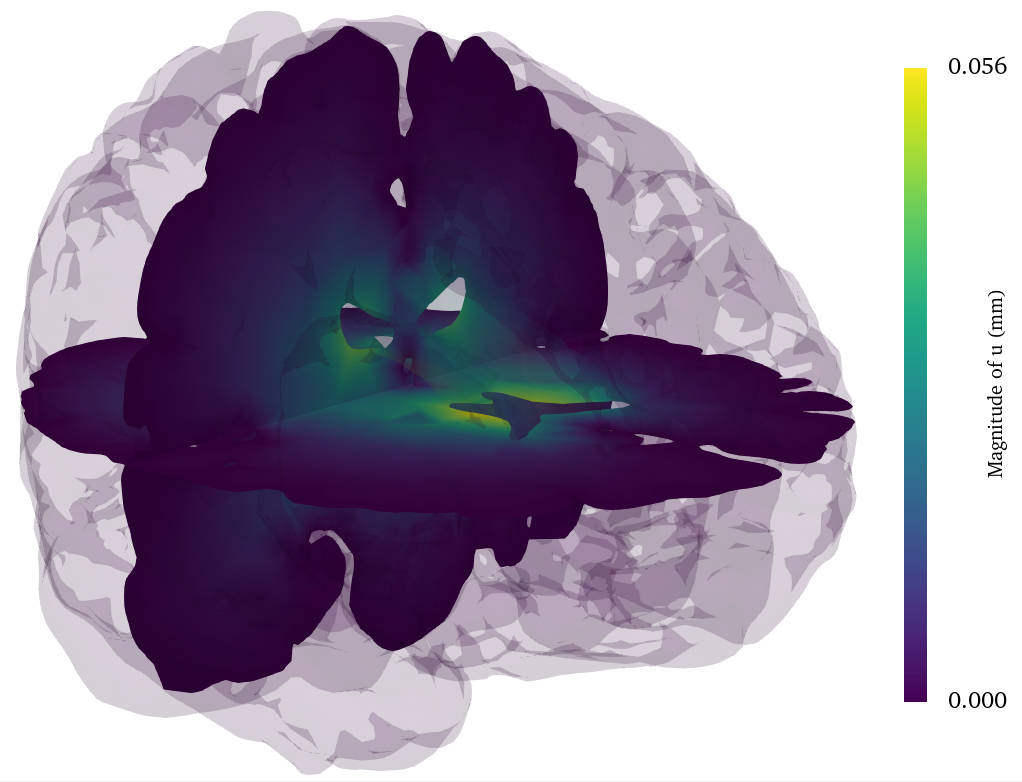}
			\caption{the coupled algorithm}
		\end{subfigure}
		\caption{Displacement magnitude values $|\bm{u}_{h}|$ computed using different algorithms.}
		\label{fig:u_distribution}
	\end{figure}

\begin{figure}[!htbp]
		\centering
		\begin{subfigure}[b]{.45\linewidth}
			\centering
			\includegraphics[width=6.0cm,height=5.0cm]{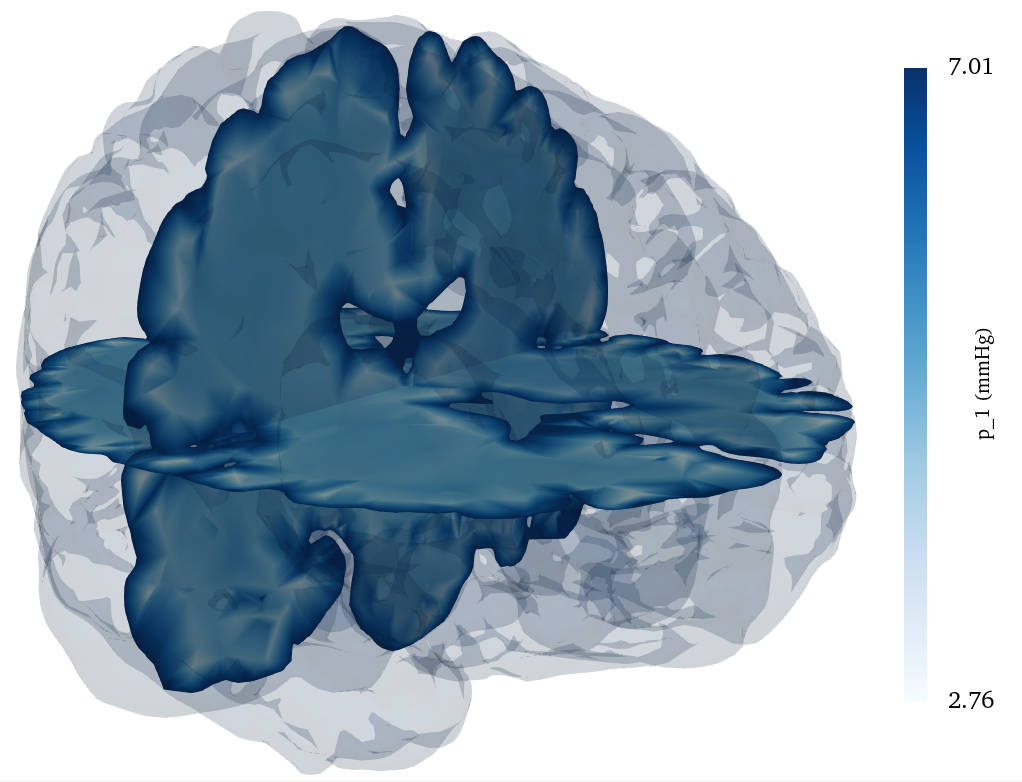}
			\caption{the iteratively decoupled algorithm}
		\end{subfigure}
        \qquad
		\begin{subfigure}[b]{.45\linewidth}
			\centering
			\includegraphics[width=6.0cm,height=5.0cm]{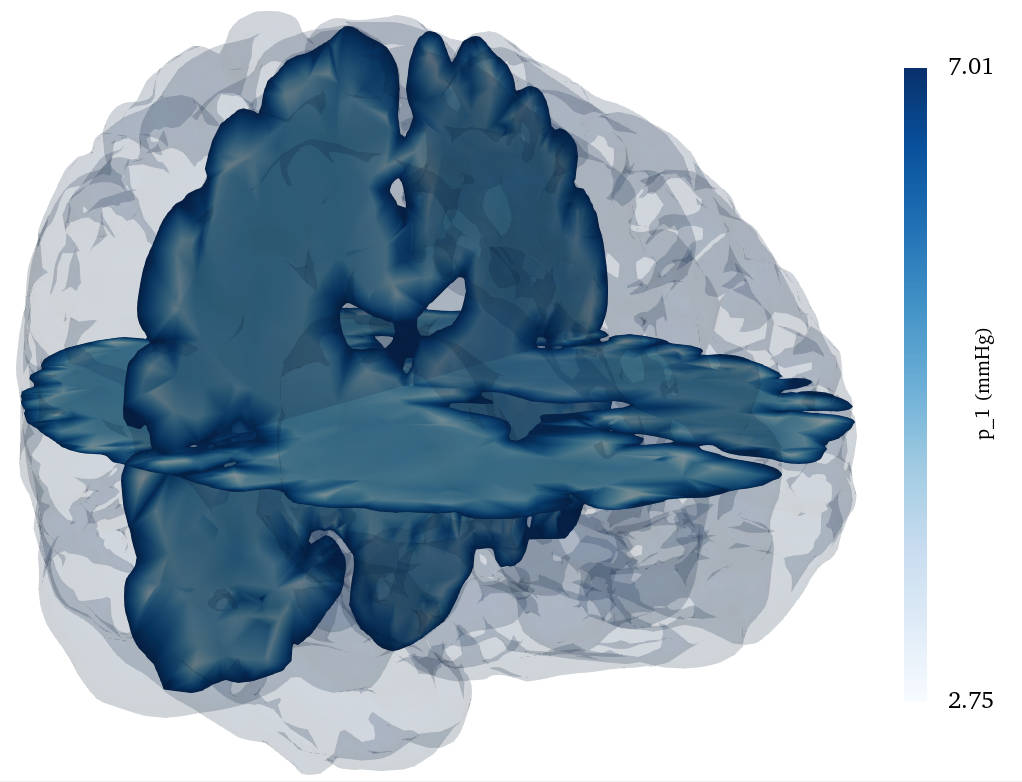}
			\caption{the coupled algorithm}
		\end{subfigure}
		\caption{Extracellular pressure values $p_{1, h}$ computed using different algorithms.}
		\label{fig:p1_distribution}
	\end{figure}

\begin{figure}[!htbp]
		\centering
		\begin{subfigure}[b]{.45\linewidth}
			\centering
			\includegraphics[width=6.0cm,height=5.0cm]{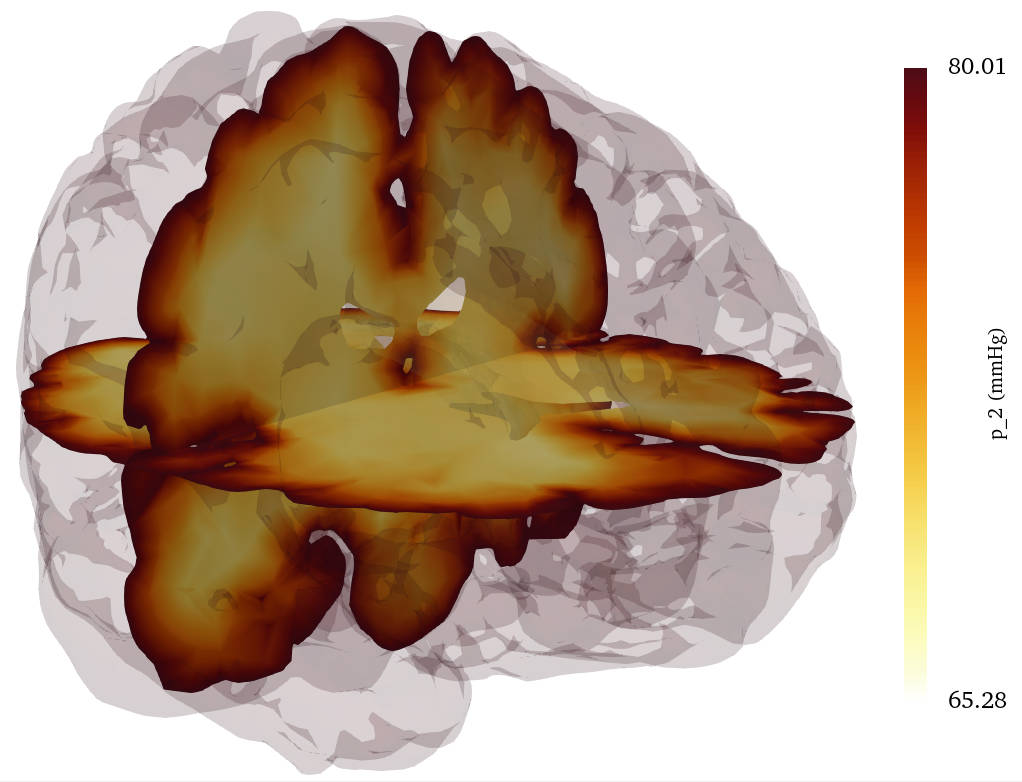}
			\caption{the iteratively decoupled algorithm}
		\end{subfigure}
        \qquad
		\begin{subfigure}[b]{.45\linewidth}
			\centering
			\includegraphics[width=6.0cm,height=5.0cm]{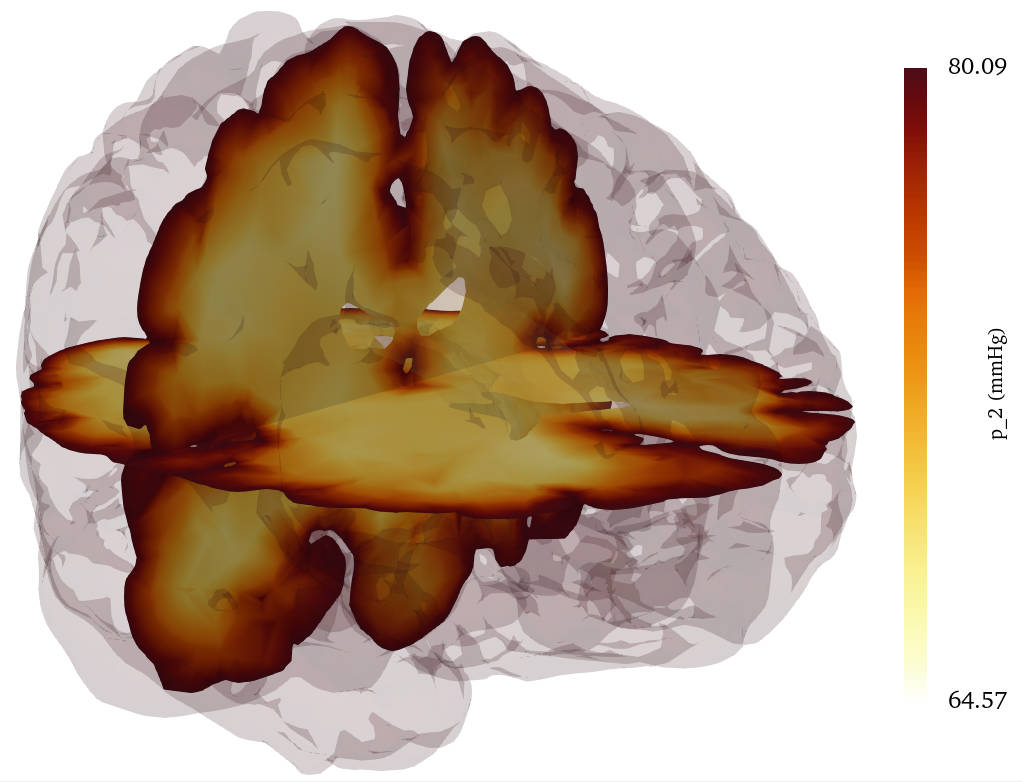}
			\caption{the coupled algorithm}
		\end{subfigure}
		\caption{Arterial pressure values $p_{2, h}$ computed using different algorithms.}
		\label{fig:p2_distribution}
	\end{figure}

\begin{figure}[!htbp]
		\centering
		\begin{subfigure}[b]{.45\linewidth}
			\centering
			\includegraphics[width=6.0cm,height=5.0cm]{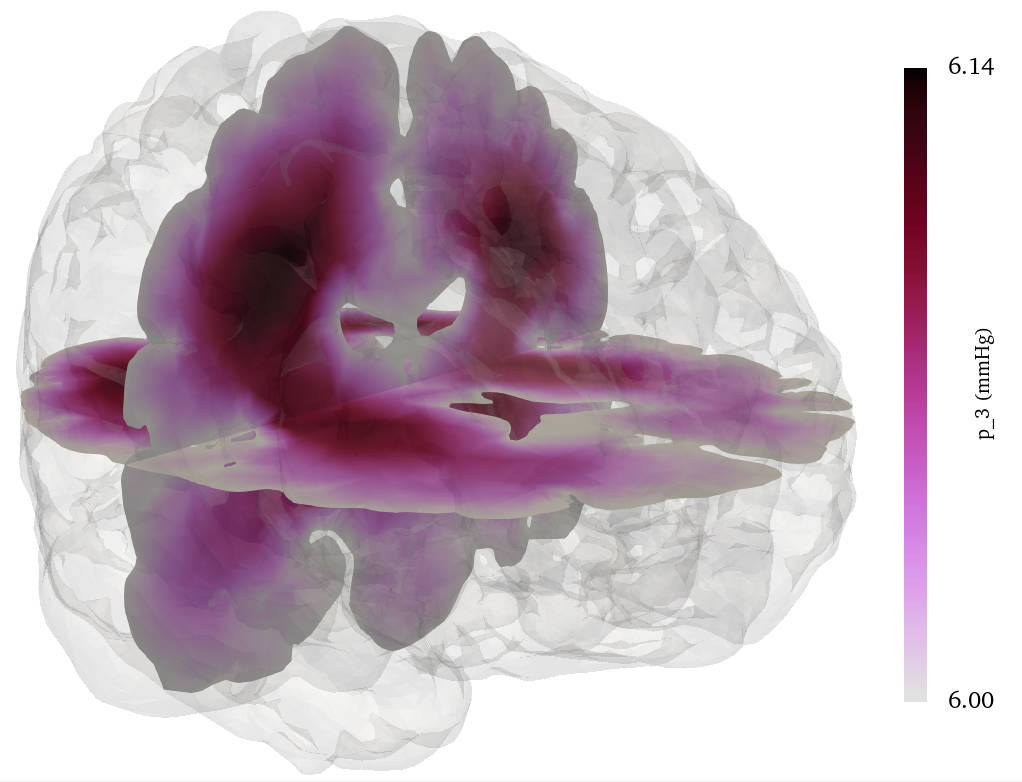}
			\caption{the iteratively decoupled algorithm}
		\end{subfigure}
        \qquad
		\begin{subfigure}[b]{.45\linewidth}
			\centering
			\includegraphics[width=6.0cm,height=5.0cm]{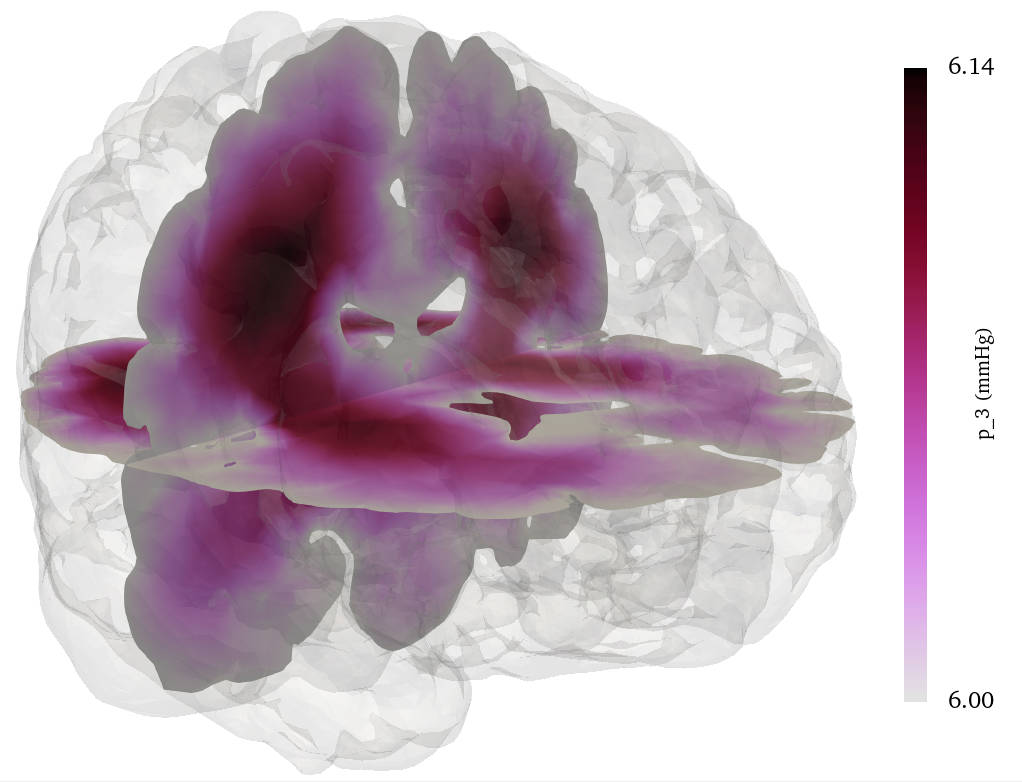}
			\caption{the coupled algorithm}
		\end{subfigure}
		\caption{Venous pressure values $p_{3, h}$ computed using different algorithms.}
		\label{fig:p3_distribution}
	\end{figure}

\begin{figure}[!htbp]
		\centering
		\begin{subfigure}[b]{.45\linewidth}
			\centering
			\includegraphics[width=6.0cm,height=5.0cm]{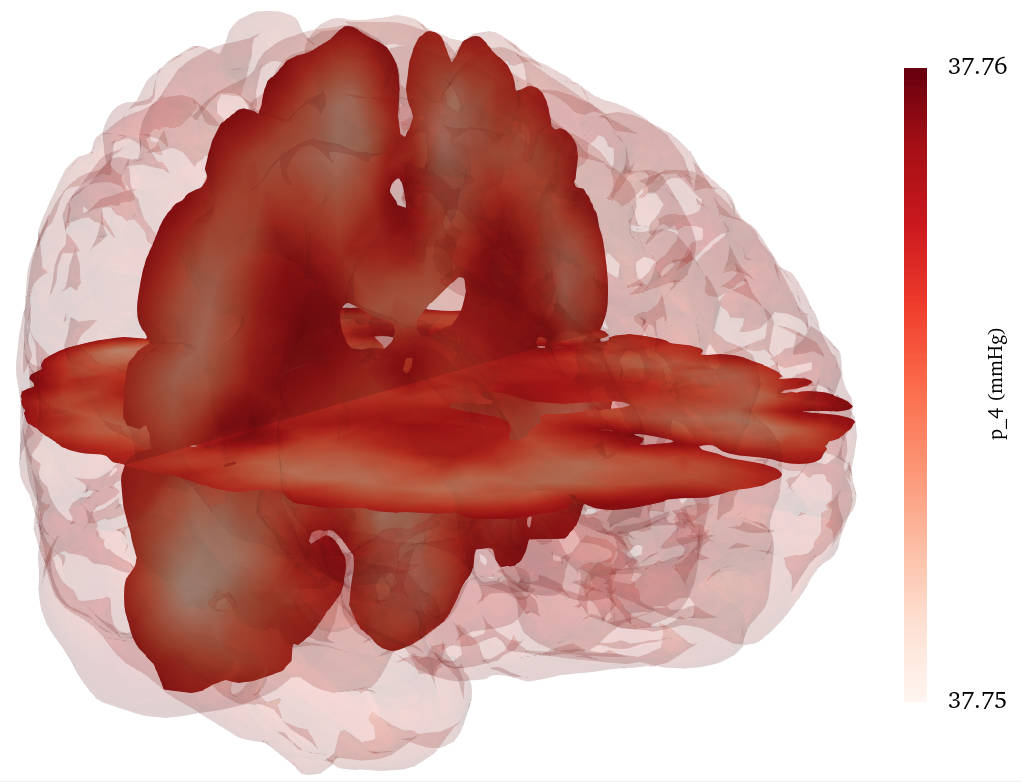}
			\caption{the iteratively decoupled algorithm}
		\end{subfigure}
        \qquad
		\begin{subfigure}[b]{.45\linewidth}
			\centering
			\includegraphics[width=6.0cm,height=5.0cm]{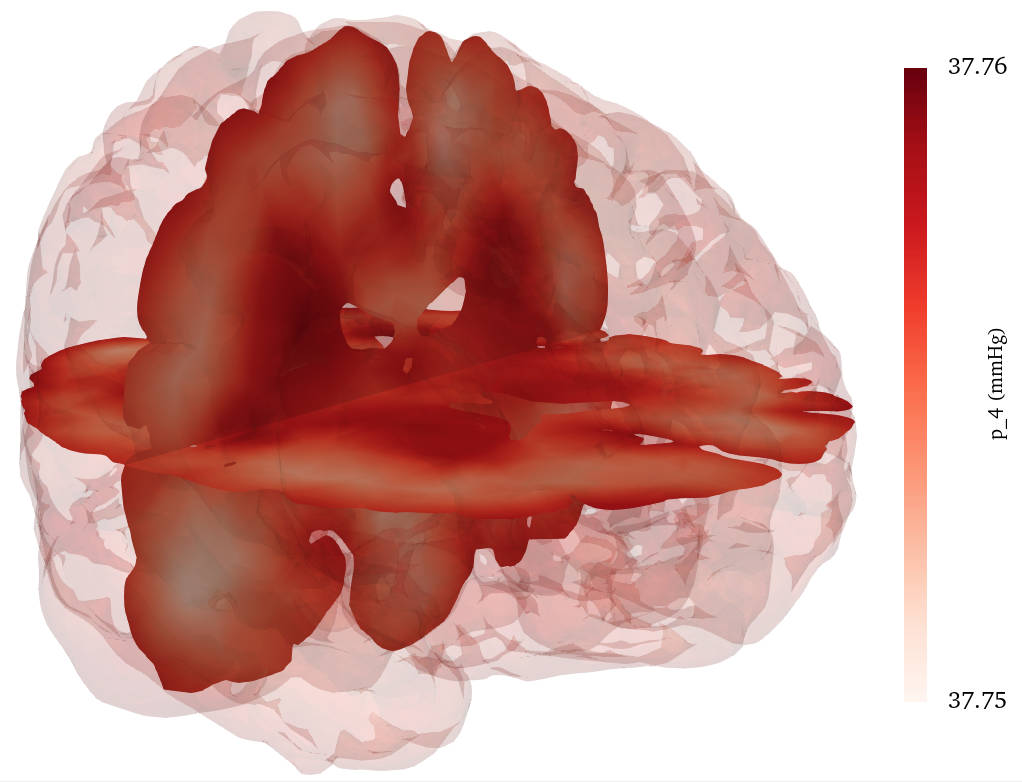}
			\caption{the coupled algorithm}
		\end{subfigure}
		\caption{Capillary pressure values $p_{4, h}$ computed using different algorithms.}
		\label{fig:p4_distribution}
	\end{figure}

	Then, we present the distribution of $|\bm{u}_{h}|$ and $p_{j, h}$, computed by 
	the iteratively decoupled algorithm and the coupled algorithm when $t = 2.25\,\textup{s}$, see Figures \ref{fig:u_distribution}--\ref{fig:p4_distribution}.

\begin{figure}[!htbp]
		\centering
		\begin{subfigure}[b]{.45\linewidth}
			\centering
			\includegraphics[width=6.0cm,height=5.0cm]{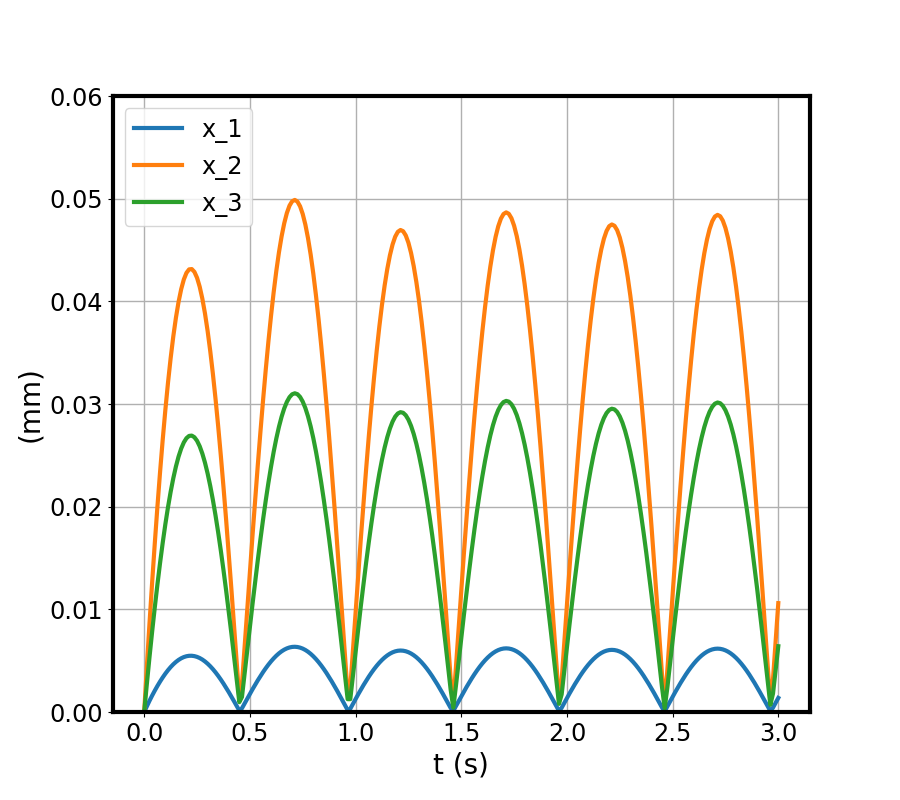}
			\caption{the iteratively decoupled algorithm}
		\end{subfigure}
		\begin{subfigure}[b]{.45\linewidth}
			\centering
			\includegraphics[width=6.0cm,height=5.0cm]{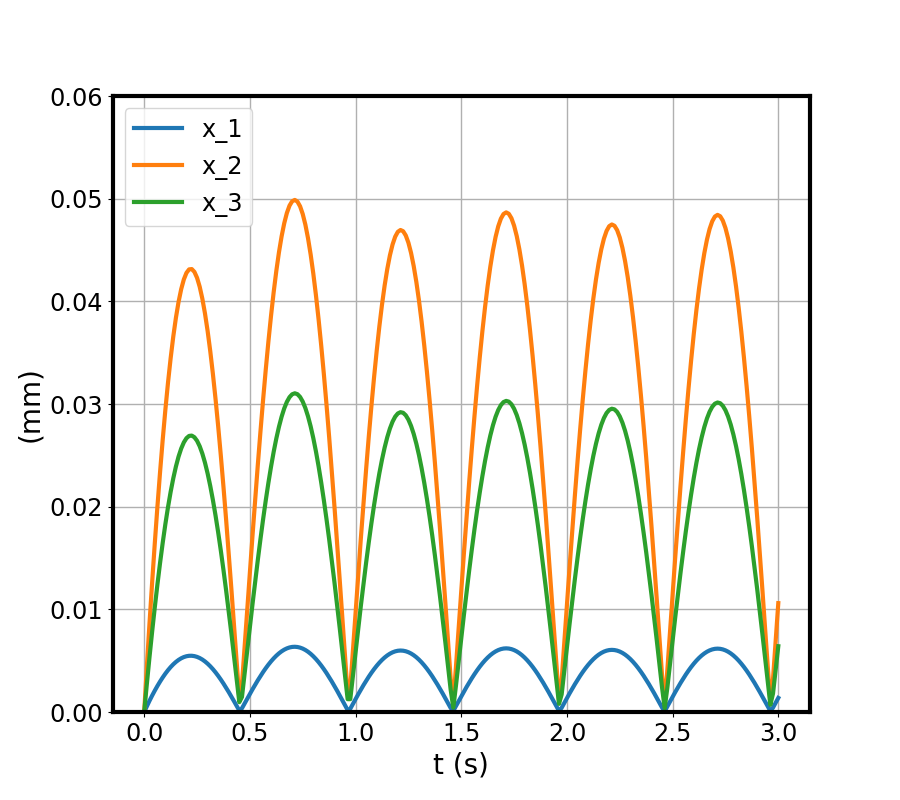}
			\caption{the coupled algorithm}
		\end{subfigure}
		\caption{Evolution of displacement magnitude values $|\bm{u}_{h}|$ at fixed points computed using different algorithms.}
		\label{fig:u_evolution}
	\end{figure}

\begin{figure}[!htbp]
		\centering
		\begin{subfigure}[b]{.45\linewidth}
			\centering
			\includegraphics[width=6.0cm,height=5.0cm]{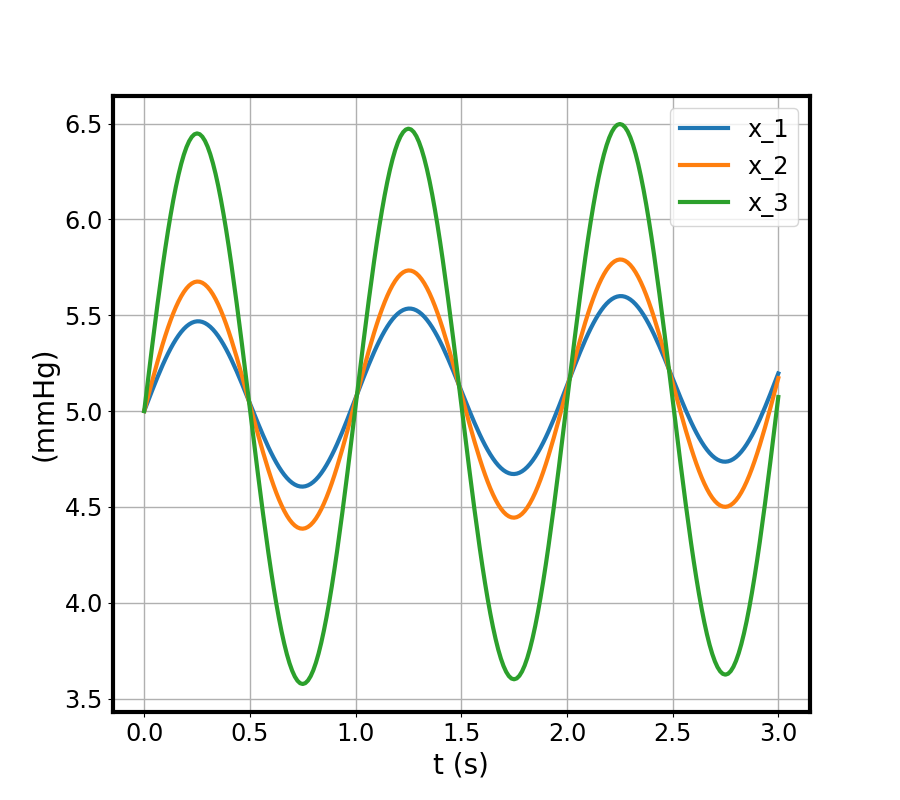}
			\caption{the iteratively decoupled algorithm}
		\end{subfigure}
		\begin{subfigure}[b]{.45\linewidth}
			\centering
			\includegraphics[width=6.0cm,height=5.0cm]{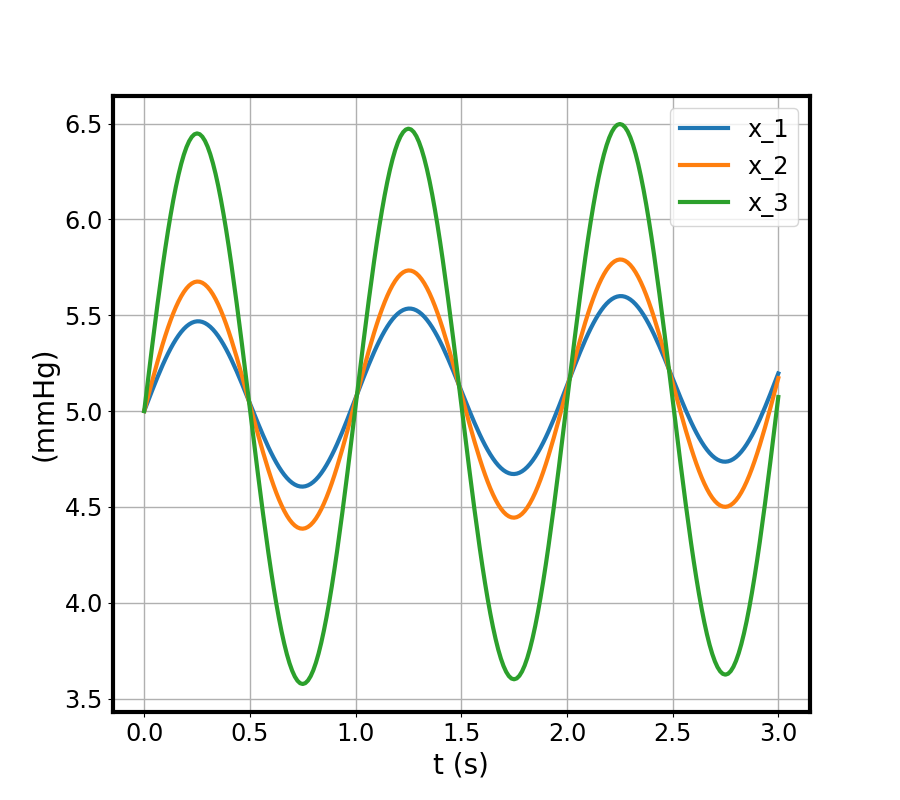}
			\caption{the coupled algorithm}
		\end{subfigure}
		\caption{Evolution of extracellular pressure values $p_{1, h}$ at fixed points computed using different algorithms.}
		\label{fig:p1_evolution}
	\end{figure}

\begin{figure}[!htbp]
		\centering
		\begin{subfigure}[b]{.45\linewidth}
			\centering
			\includegraphics[width=6.0cm,height=5.0cm]{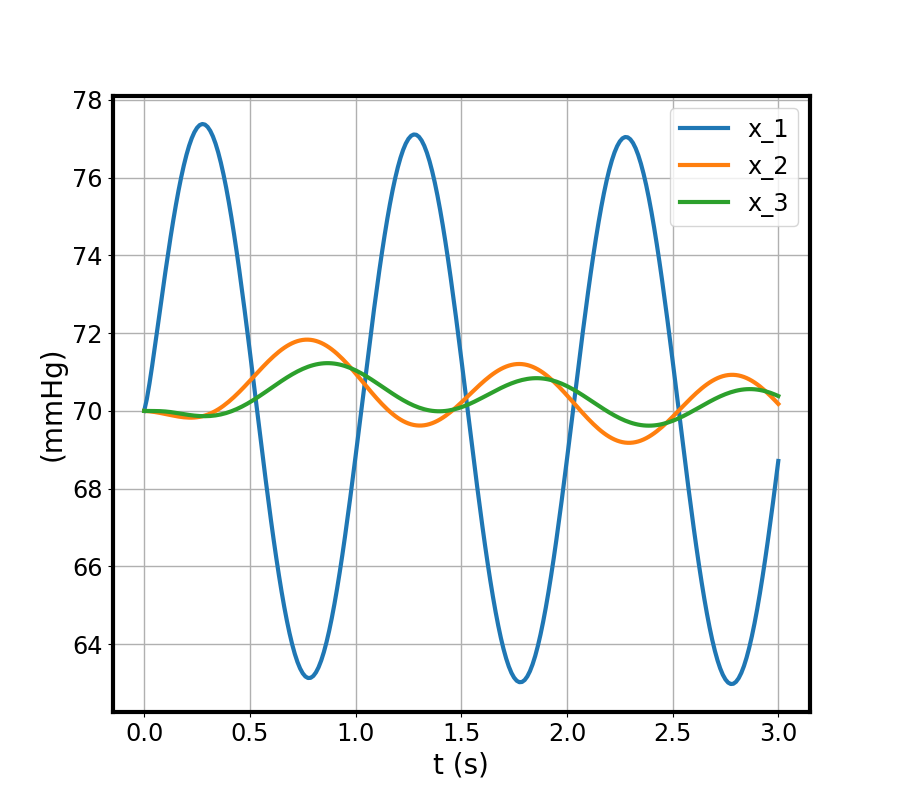}
			\caption{the iteratively decoupled algorithm}
		\end{subfigure}
		\begin{subfigure}[b]{.45\linewidth}
			\centering
			\includegraphics[width=6.0cm,height=5.0cm]{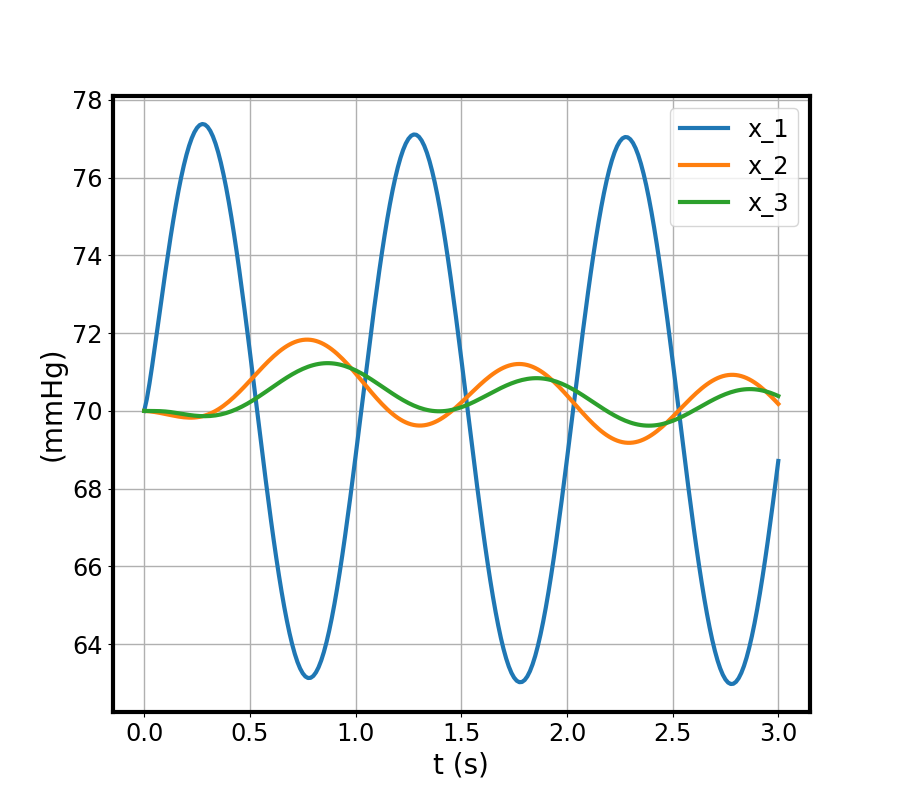}
			\caption{the coupled algorithm}
		\end{subfigure}
		\caption{Evolution of arterial pressure values $p_{2, h}$ at fixed points computed using different algorithms.}
		\label{fig:p2_evolution}
	\end{figure}

\begin{figure}[!htbp]
		\centering
		\begin{subfigure}[b]{.45\linewidth}
			\centering
			\includegraphics[width=6.0cm,height=5.0cm]{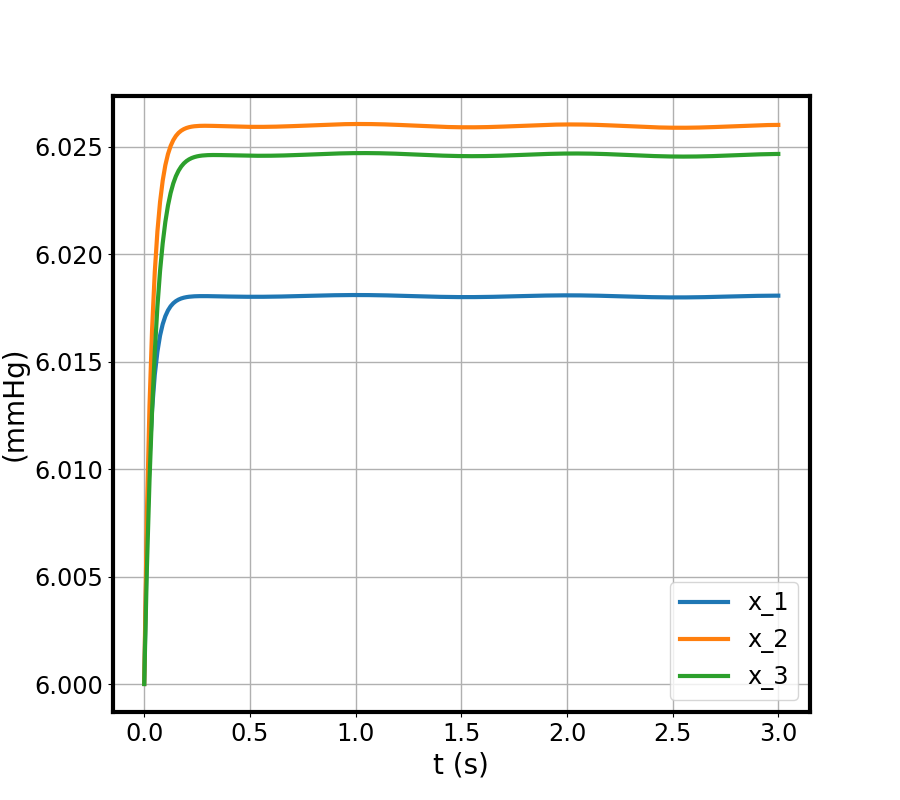}
			\caption{the iteratively decoupled algorithm}
		\end{subfigure}
		\begin{subfigure}[b]{.45\linewidth}
			\centering
			\includegraphics[width=6.0cm,height=5.0cm]{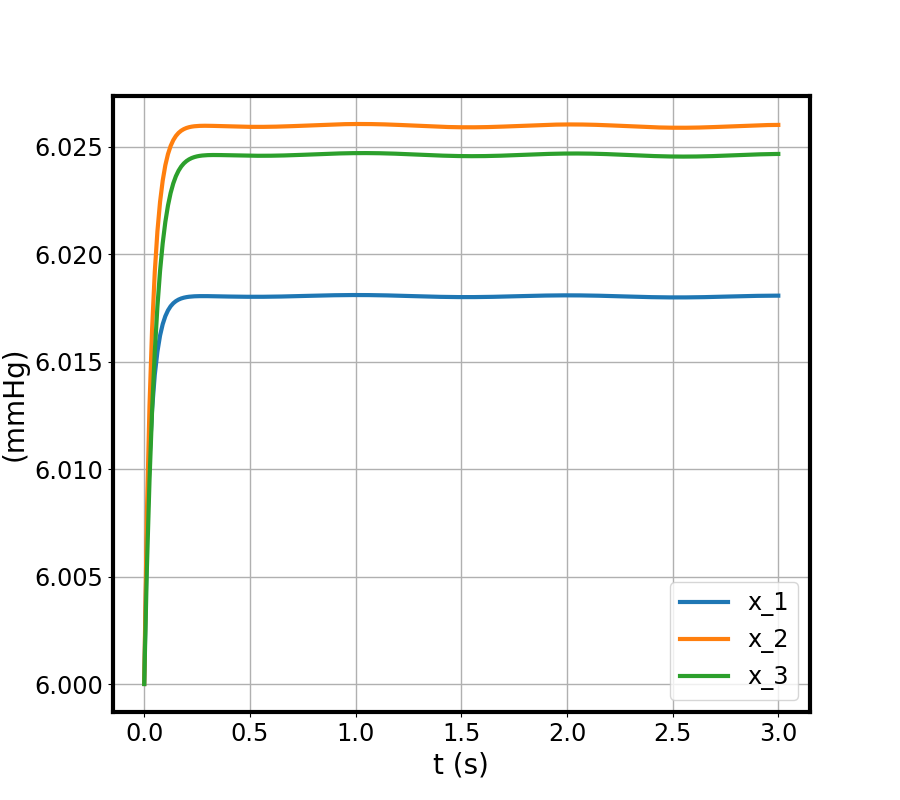}
			\caption{the coupled algorithm}
		\end{subfigure}
		\caption{Evolution of venous pressure values $p_{3, h}$ at fixed points computed using different algorithms.}
		\label{fig:p3_evolution}
	\end{figure}

\begin{figure}[!htbp]
		\centering
		\begin{subfigure}[b]{.45\linewidth}
			\centering
			\includegraphics[width=6.0cm,height=5.0cm]{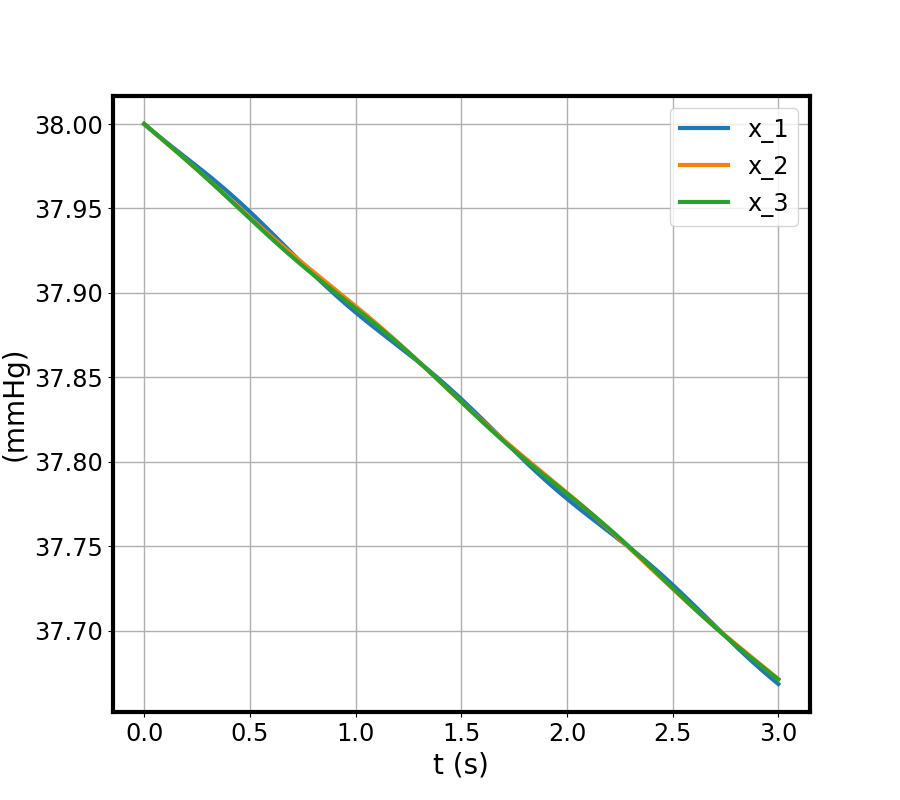}
			\caption{the iteratively decoupled algorithm}
		\end{subfigure}
		\begin{subfigure}[b]{.45\linewidth}
			\centering
			\includegraphics[width=6.0cm,height=5.0cm]{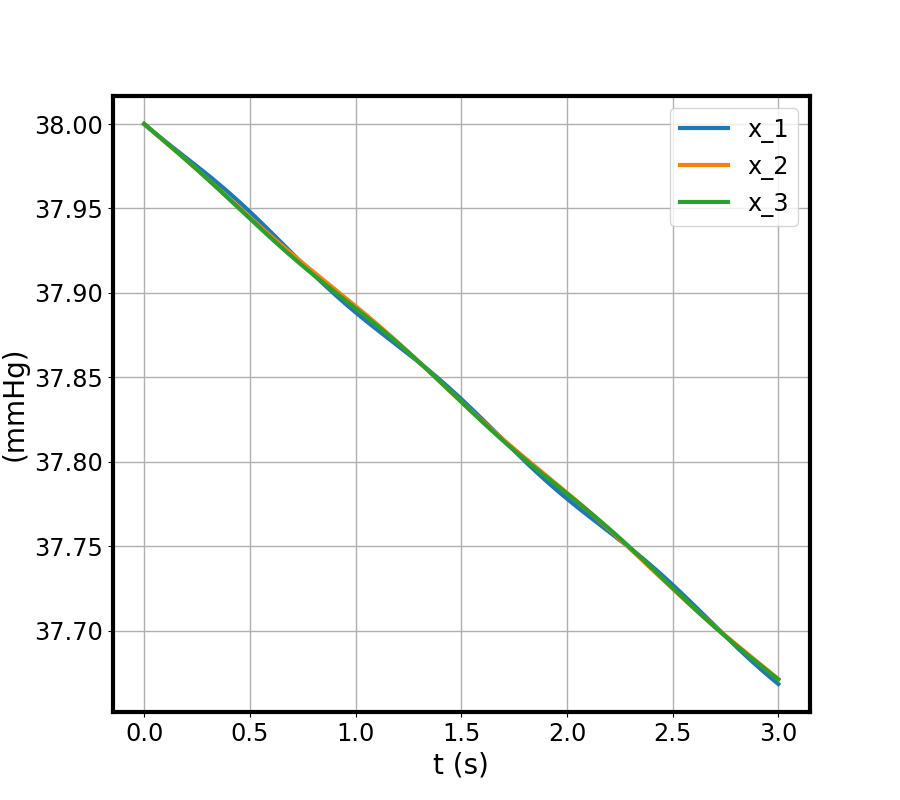}
			\caption{the coupled algorithm}
		\end{subfigure}
		\caption{Evolution of capillary pressure values $p_{4, h}$ at fixed points computed using different algorithms.}
		\label{fig:p4_evolution}
	\end{figure}

	We also calculate the time evolution of the numerical solution values at three fixed points 
	$\bm{x}_{1} = (89.9, 108.9, 82.3)$ (center), $\bm{x}_{2} = (102.2, 139.3, 82.3)$ (point in the central $z$-plane), 
	and $\bm{x}_{3} = (110.7, 108.9, 98.5)$ (point in the central $y$-plane), see Figures \ref{fig:u_evolution}--\ref{fig:p4_evolution}. 
    For ease of comparison, the same time step size is used for both algorithms, i.e., $\Delta t_{1} = \Delta t_{2} = 0.0125\,\textup{s}$.

	As shown in Figures \ref{fig:u_distribution}--\ref{fig:p4_evolution}, 
	the numerical results computed by the iteratively decoupled algorithm closely match those of the coupled algorithm 
    as well as the results presented in Section 6 of \cite{lee2019mixed}, validating the effectiveness and efficiency of the iteratively decoupled algorithm.


\section{Conclusions}
In this work, we have developed an iteratively decoupled algorithm for solving multiple-network poroelasticity problems. Our theoretical analysis demonstrates that the sequence ${(\pmb{u}_h^{n,k},\xi_h^{n,k},\vec{p}_h^{n,k})}$ generated by the decoupled scheme converges to the coupled solution ${(\pmb{u}_h^n,\xi_h^n,\vec{p}_h^n)}$, with guaranteed unconditional stability and convergence. 
Numerical experiments confirm the validity of our theoretical results and demonstrate the algorithm's computational efficiency. When applied to brain edema simulations, the proposed method achieves comparable accuracy to the fully coupled approach while significantly reducing computational costs. These findings suggest that our decoupled algorithm offers an effective balance between accuracy and efficiency for complex multi-network poroelastic systems.  



\section*{Acknowledgments}
The work of M. Cai is partially supported by the NIH-RCMI award (Grant No. 347U54MD013376), the Army Research Office award (W911NF2310004), and the affiliated project award from the Center for Equitable Artificial Intelligence and Machine Learning Systems (CEAMLS) at Morgan State University (project ID 02232301). 
The work of M. Lei and F. Wang is partially supported by NSFC (Grant No. 12071227, 12371369), and the Ministry of Science and Technology of China (Grant No. 2020YFA0713800). The work of J. Li and J. Sun is partially supported by the Shenzhen Sci-Tech Fund No. RCJC20200714114556020, Guangdong Basic and Applied Research Fund No. 2023B1515250005,  National Center for Applied Mathematics Shenzhen, and SUSTech International Center for Mathematics.

\section*{Data Availability}
The code for the numerical tests is available upon request.

\section*{Declarations}
\noindent \textbf{Conflict of interest}
The authors have no competing interests.

\begin{appendices}
\section{Error analysis of the coupled algorithm}
\setcounter{section}{0}
\setcounter{equation}{0}
\setcounter{theorem}{0}

\renewcommand{\theequation}{A.\arabic{equation}}
\renewcommand{\thesection}{A.\arabic{section}}
\renewcommand{\thetheorem}{A.\arabic{theorem}}
\renewcommand{\thelemma}{A.\arabic{lemma}}

This appendix provides the optimal order error estimate for the coupled algorithm. Specifically, we demonstrate that the temporal error is of order $O(\Delta t)$; The energy-norm errors for $\pmb{u}$ and $\vec{p}$ are of order $O(h^2)$ and order $O(h)$, respectively;
The $L^2$ norm errors for $\xi$ and $\vec{p}$ are of order $O(h^2)$. 
On the continuous level, we refer readers to \cite{lee2019mixed} for the energy estimate of the variational problem \eqref{weak_problem}. In the following lemma, we provide the energy estimate for the coupled algorithm.
\vspace{0.5cm}

\begin{lemma}
   Let $\{(\pmb{u}_h^n,\xi_h^n,\vec{p}_h^n) \}$ be the solutions of problem \eqref{coupled_alg1}-\eqref{coupled_alg3}, then the following identity holds
   \begin{align} \label{coupled_alg_conv1}
    J_h^l+S_h^l=J_h^0\quad for\,l \in\{1,2,\cdots,L\}.
\end{align}
where
\begin{align*}
  & J_h^l:=\mu |\!|\epsilon(\pmb{u}_h^l)|\!|_{L^2(\Omega)}^2+\frac{1}{2\lambda}|\!|\vec{\alpha}^{\mathsf{T}} \vec{p}_h^l-\xi_h^l|\!|_{L^2(\Omega)}^2+\frac{1}{2}|\!|\hat{S}\vec{p}_h^l|\!|_{L^2(\Omega)}^2-(\pmb{f},\pmb{u}_h^l)-\langle \pmb{h},\pmb{u}_h^l\rangle_{\Gamma_{\pmb{u},N}}. \\
 &  S_h^l:= \\ \notag
 &\Delta t \sum_{n=1}^l\bigg[\Delta t \left( \mu |\!|d_t \epsilon(\pmb{u}_h^n)|\!|_{L^2(\Omega)}^2+\frac{1}{2\lambda}|\!|d_t(\vec{\alpha}^{\mathsf{T}} \vec{p}_h^n-\xi_h^n)|\!|_{L^2(\Omega)}^2+\frac{1}{2}|\!|d_t (\hat{S}\vec{p}_h^n)|\!|_{L^2(\Omega)}^2 \right)\\ \notag
  &+|\!|\hat{K}\nabla \vec{p}_h^n|\!|_{L^2(\Omega)}^2+\sum_{1\leq i<j\leq N}\beta_{i,j} |\!|p_{i,h}^n-p_{j,h}^n|\!|_{L^2(\Omega)}^2-(\vec{g},\vec{p}_h^n)-\left\langle \vec{l},\vec{p}_h^n\right\rangle_{\Gamma_{\vec{p},N}} \bigg].
\end{align*}
Here, we denote $d_t \eta^n :=\frac{\eta^n-\eta^{n-1}}{\Delta t}$, where $\eta$ can be a vector or a scalar, $\hat{S}=\text{diag}(\sqrt{c_1},\sqrt{c_2},\cdots,\sqrt{c_N}), \hat{K}=\text{diag}(\sqrt{K_1},\sqrt{K_2},\cdots,\sqrt{K_N})$. Moreover, there holds
\begin{align} \label{appendix_lemma_con1}
    |\!|\xi_h^l|\!|_{L^2(\Omega)} \leq C (|\!|\epsilon(\pmb{u}_h^l)|\!|_{L^2(\Omega)}+|\!|\pmb{f}|\!|_{L^2(\Omega)}+|\!|\pmb{h}|\!|_{L^2(\Gamma_{\pmb{u},N})} ),
\end{align}
where $C$ is a positive constant.
\end{lemma}

\noindent{\bf Proof}. Setting $\pmb{v}_h=d_t \pmb{u}_h^n$ in \eqref{coupled_alg1}, $\eta_h=\xi_h^n$ in \eqref{coupled_alg2}, $\vec{q}_{h}=\vec{p}_{h}^n$ in \eqref{coupled_alg3}, we have
\begin{align}
\label{coupled_alg_conv2}
    2\mu (\epsilon(\pmb{u}_h^n),d_t \epsilon(\pmb{u}_h^n))-(\xi_h^n,\operatorname{div} d_t \pmb{u}_h^n)&=d_t (\pmb{f},\pmb{u}_h^n)+d_t \langle \pmb{h},\pmb{u}_h^n\rangle_{\Gamma_{\pmb{u},N}},\\
\label{coupled_alg_conv3}
    (\operatorname{div} d_t \pmb{u}_h^n,\xi_h^n)+\frac{1}{\lambda}(d_t \xi_h^n,\xi_h^n)&=\frac{1}{\lambda}(\vec{\alpha}^{\mathsf{T}} d_t \vec{p}_h^n,\xi_h^n),\\
\label{coupled_alg_conv4}
    \left( \left( S+\frac{1}{\lambda}\vec{\alpha}\vec{\alpha}^{\mathsf{T}} \right)d_t \vec{p}_h^n ,\vec{p}_h^n \right)-\frac{1}{\lambda}\left( \vec{\alpha}d_t \xi_h^n ,\vec{p}_h^n    \right)&+\left(K\nabla \vec{p}_h^n, \nabla \vec{p}_h^n\right)+\left(B\vec{p}_h^n,\vec{p}_h^n\right) \\ \notag
    &=(\vec{g},\vec{p}_h^n)+\left\langle \vec{l},\vec{p}_h^n\right\rangle_{\Gamma_{\vec{p},N}}.
\end{align}
Summing up \eqref{coupled_alg_conv2}-\eqref{coupled_alg_conv4}, we get
\begin{align} \label{coupled_alg_conv5}
  & 2\mu(\epsilon(\pmb{u}_h^n),d_t \epsilon(\pmb{u}_h^n))-d_t (\pmb{f},\pmb{u}_h^n)-d_t \langle \pmb{h},\pmb{u}_h^n\rangle_{\Gamma_{\pmb{u},N}}+\frac{1}{\lambda}(d_t \xi_h^n,\xi_h^n)-\frac{1}{\lambda}(\vec{\alpha}^{\mathsf{T}} d_t \vec{p}_h^n,\xi_h^n)\\ \notag
 &+\left( \left( S+\frac{1}{\lambda}\vec{\alpha}\vec{\alpha}^{\mathsf{T}} \right)d_t \vec{p}_h^n ,\vec{p}_h^n \right)-\frac{1}{\lambda}\left( \vec{\alpha}d_t \xi_h^n ,\vec{p}_h^n    \right)+\left(K\nabla \vec{p}_h^n, \nabla \vec{p}_h^n\right)+\left(B\vec{p}_h^n,\vec{p}_h^n\right)-(\vec{g},\vec{p}_h^n)\\ \notag
 &-\left\langle \vec{l},\vec{p}_h^n\right\rangle_{\Gamma_{\vec{p},N}}=0.
\end{align}
Note that, the following identities and  \eqref{identity_1} hold.
\begin{align*}
    2(\eta_h^n,d_t \eta_h^n)=d_t |\!|\eta_h^n|\!|_{L^2(\Omega)}^2&+\Delta t |\!|d_t \eta_h^n|\!|_{L^2(\Omega)}^2,\\
   \frac{1}{\lambda}(d_t \xi_h^n,\xi_h^n)-\frac{1}{\lambda}(\vec{\alpha}^{\mathsf{T}} d_t \vec{p}_h^n,\xi_h^n)&+\frac{1}{\lambda}\left( \vec{\alpha}\vec{\alpha}^{\mathsf{T}} d_t \vec{p}_h^n,\vec{p}_h^n  \right)-\frac{1}{\lambda}\left( \vec{\alpha}d_t\xi_h^n,\vec{p}_h^n \right) \\ \notag
 &=\frac{1}{\lambda}(\vec{\alpha}^{\mathsf{T}}\vec{p}_h^n-\xi_h^n,d_t(\vec{\alpha}^{\mathsf{T}} \vec{p}_h^n-\xi_h^n)).
  \end{align*}
According to the above identities, \eqref{coupled_alg_conv5} can be rewritten as
\begin{align} \label{coupled_alg_conv6}
   & d_t \bigg( \mu |\!|\epsilon(\pmb{u}_h^n)|\!|_{L^2(\Omega)}^2+\frac{1}{2\lambda}|\!|\vec{\alpha}^{\mathsf{T}} \vec{p}_h^n-\xi_h^n|\!|_{L^2(\Omega)}^2+\frac{1}{2}|\!|\hat{S}\vec{p}_h^n|\!|_{L^2(\Omega)}^2
    -(\pmb{f},\pmb{u}_h^n) \\ \notag
    &-\langle \pmb{h},\pmb{u}_h^n\rangle_{\Gamma_{\pmb{u},N}} \bigg)
    +\Delta t \left( \mu |\!|d_t \epsilon(\pmb{u}_h^n)|\!|_{L^2(\Omega)}^2+\frac{1}{2\lambda}|\!|d_t(\vec{\alpha}^{\mathsf{T}} \vec{p}_h^n-\xi_h^n)|\!|_{L^2(\Omega)}^2+\frac{1}{2}|\!|d_t (\hat{S}\vec{p}_h^n)|\!|_{L^2(\Omega)}^2 \right) \\ \notag
 & +|\!|\hat{K}\nabla \vec{p}_h^n|\!|_{L^2(\Omega)}^2
+\sum_{1\leq i<j\leq N}\beta_{i,j} |\!|p_{i,h}^n-p_{j,h}^n|\!|_{L^2(\Omega)}^2-(\vec{g},\vec{p}_h^n)-\left\langle \vec{l},\vec{p}_h^n\right\rangle_{\Gamma_{\vec{p},N}}=0.
\end{align}
Applying the summation operator $\Delta t \sum_{n=1}^l$ to both sides of the \eqref{coupled_alg_conv6}, we get \eqref{coupled_alg_conv1}. In addition, using \eqref{discrete inf-sup condition} and \eqref{coupled_alg1}, we have \eqref{appendix_lemma_con1}.
$\hfill\qedsymbol$

\vspace{0.5cm}
We introduce a set of auxiliary projection operators as a preliminary step for a priori estimate of the discrete formulation \eqref{coupled_alg}.  In particular, we define projection operators
$$\Pi_h^{\pmb{V}}:\pmb{V}\rightarrow \pmb{V}_h,~~~~\Pi_h^W:W\rightarrow W_h,$$
as follows.
For any $(\pmb{u},\xi)\in \pmb{V}\times W$, we define its projection $(\Pi_h^{\pmb{V}} \pmb{u},\Pi_h^W \xi)\in \pmb{V}_h \times W_h$ as the unique discrete solution to the Stokes problem:
\begin{subequations}
\label{Stokes_projection}
\begin{align}
\label{Stokes_projection_eq1}
   2\mu(\epsilon(\Pi_h^{\pmb{V}} \pmb{u}),\epsilon(\pmb{v}_h))-(\Pi_h^W \xi,\operatorname{div}\pmb{v}_h)&=2\mu(\epsilon(\pmb{u}),\epsilon(\pmb{v}_h))-(\xi,\operatorname{div}\pmb{v}_h) \quad \forall \pmb{v}_h \in \pmb{V}_h,\\
 \label{Stokes_projection_eq2}
-(\operatorname{div}\Pi_h^{\pmb{V}}\pmb{u},\eta_h)&=-(\operatorname{div}\pmb{u},\eta_h)\quad
\forall \eta_h \in  W_h.
\end{align}
\end{subequations}
We also define the projection operator
$$\Pi_h^{M} :M\rightarrow M_h,$$
as follows.
For any $\vec{p} \in M$, we define its projection $\Pi_h^{M} \vec{p}\in M_{h}$ as the unique discrete solution to the weighted elliptic problem
\begin{align} \label{ell_projection_op}
    (K \nabla \Pi_h^{M} \vec{p},\nabla \vec{q}_{h})= (K \nabla \vec{p},\nabla \vec{q}_{h})\quad \forall \vec{q}_h \in M_h.
\end{align}
Here, we list the properties of the operators $(\Pi_h^{\pmb{V}},\Pi_h^W,\Pi_h^{M})$. For all $(\pmb{u},\xi,\vec{p})\in [H^3(\Omega)]^2\times H^2(\Omega)\times [H^2(\Omega)]^N$, there holds \cite{lee2019mixed}
\begin{align} \label{projection_err1}
    |\!|\pmb{u}-\Pi_h^{\pmb{V}}  \pmb{u}|\!|_{H^1 (\Omega)}+|\!|\xi-\Pi_h^W \xi|\!|_{L^2 (\Omega)}&\leq Ch^2 \left( |\!|\pmb{u}|\!|_{H^3 (\Omega)}+|\!|\xi|\!|_{H^2 (\Omega)}\right),\\
\label{projection_err2}
    |\!|\vec{p} -\Pi_h^{M} \vec{p} |\!|_{H^1 (\Omega)}&\leq Ch|\!|\vec{p}|\!|_{H^2(\Omega)},\\
\label{projection_err3}
    |\!|\vec{p} -\Pi_h^{M} \vec{p} |\!|_{L^2 (\Omega)}&\leq Ch^2|\!|\vec{p}|\!|_{H^2(\Omega)}.
\end{align}

For convenience,  we introduce the standard decomposition of the errors into projection and discretization errors in the subsequent analysis.
\begin{align*}
    e_{\pmb{u}}^n=\pmb{u}^n-\pmb{u}_h^n=(\pmb{u}^n-\Pi_h^{\pmb{V}} \pmb{u}^n)+(\Pi_h^{\pmb{V}} \pmb{u}^n-\pmb{u}_h^n):=e_{\pmb{u}}^{I,n}+e_{\pmb{u}}^{h,n},\\
    e_{\xi}^n=\xi^n-\xi_h^n=(\xi^n-\Pi_h^W \xi^n)+(\Pi_h^W \xi^n-\xi_h^n):=e_{\xi}^{I,n}+e_{\xi}^{h,n},\\
     e_{\vec{p}}^n=\vec{p}^n-\vec{p}_{h}^n=(\vec{p}^n-\Pi_h^{M}\vec{p}^n)+(\Pi_h^{M} \vec{p}^n-\vec{p}_{h}^n):=e_{\vec{p}}^{I,n}+e_{\vec{p}}^{h,n}.
\end{align*}

\begin{lemma}
     Let $\{(\pmb{u}_h^n,\xi_h^n,\vec{p}_h^n) \}$ be the solutions of problem \eqref{coupled_alg1}-\eqref{coupled_alg3}, then the following identity holds
\begin{align} \label{coupled_alg_conv7}
    & E_h^l  +\Delta t\sum_{n=1}^l \left( |\!|\hat{K}\nabla e_{\vec{p}}^{h,n}|\!|_{L^2(\Omega)}^2+\left( Be_{\vec{p}}^{h,n},e_{\vec{p}}^{h,n}\right)\right) \\ \notag
   &+(\Delta t)^2 \sum_{n=1}^l\left(\mu|\!|d_t \epsilon(e_{\pmb{u}}^{h,n})|\!|_{L^2(\Omega)}^2+\frac{1}{2\lambda}|\!|d_t(\vec{\alpha}^{\mathsf{T}} e_{\vec{p}}^{h,n}-e_{\xi}^{h,n})|\!|_{L^2(\Omega)}^2+\frac{1}{2}|\!|d_t \hat{S}e_{\vec{p}}^{h,n}|\!|_{L^2(\Omega)}^2\right)\\ \notag
 &=E_h^0+\Delta t \sum_{n=1}^l \Bigg[ \left(\operatorname{div}(d_t \pmb{u}^n-\dot{\pmb{u}}^n),e_{\xi}^{h,n}\right)
  +\frac{1}{\lambda}\left(d_t \Pi_h^W \xi^n-\dot{\xi}^n,e_{\xi}^{h,n}\right) \\ \notag
  &-\frac{1}{\lambda}\left(\vec{\alpha}^{\mathsf{T}} (d_t \Pi_h^{M} \vec{p}^n-\dot{\vec{p}}^n),e_{\xi}^{h,n}\right)
  +\left(     B\left( \Pi_h^{M} \vec{p}^n-\vec{p}^n \right) ,e_{\vec{p}}^{h,n}\right)\\ \notag
&-\frac{1}{\lambda}\left( \vec{\alpha}\left( d_t \Pi_h^W \xi^n-\dot{\xi}^n \right)  , e_{\vec{p}}^{h,n} \right)
 +\bigg(  \bigg(  S+\frac{1}{\lambda}\vec{\alpha}\vec{\alpha}^{\mathsf{T}} \bigg) \left(d_t \Pi_h^{M}\vec{p}^n-\dot{\vec{p}}^n   \right), e_{\vec{p}}^{h,n}  \bigg) \Bigg],
\end{align}
where
$$E_h^l :=\mu|\!|\epsilon(e_{\pmb{u}}^{h,l})|\!|_{L^2(\Omega)}^2+\frac{1}{2\lambda}|\!|\vec{\alpha} ^{\mathsf{T}} e_{\vec{p}}^{h,l}-e_{\xi}^{h,l}|\!|_{L^2(\Omega)}^2+\frac{1}{2}|\!|\hat{S}e_{\vec{p}}^{h,l}|\!|_{L^2(\Omega)}^2.$$
\end{lemma}

\noindent{\bf Proof}. Letting $t=t_n,\pmb{v}=\pmb{v}_h$ in \eqref{three_field_form_wp1}, subtracting \eqref{coupled_alg1} from equation \eqref{three_field_form_wp1} and combining with  \eqref{Stokes_projection}, there holds
\begin{align} \label{coupled_alg_conv8}
    2 \mu (\epsilon(e_{\pmb{u}}^{h,n}),\epsilon(\pmb{v}_h))-(e_{\xi}^{h,n},\operatorname{div}\pmb{v}_h)=0,
\end{align}

Taking the derivative of \eqref{three_field_form_wp2}  with respect to time $t$, and letting $t=t_n,\eta=\eta_h$, we have
\begin{align} \label{lemma5.7_eq1}
     (\text{div}\dot{\pmb{u}}^n,\eta_h)+\frac{1}{\lambda}(\dot{\xi}^n,\eta_h)-\frac{1}{\lambda}(\vec{\alpha}^{\mathsf{T}} \dot{\vec{p}}^n,\eta_h)=0.
\end{align}
Using \eqref{coupled_alg2}, it holds that
\begin{align} \label{lemma5.7_eq2}
    (\operatorname{div} d_t \pmb{u}^n,\eta_h)+\frac{1}{\lambda}(d_t \xi^n,\eta_h)-\frac{1}{\lambda}(\vec{\alpha}^{\mathsf{T}} d_t \vec{p}^n,\eta_h)=0.
\end{align}
By subtracting \eqref{lemma5.7_eq1} from \eqref{lemma5.7_eq2}, we obtain
\begin{align} \label{coupled_alg_conv9}
    (\operatorname{div}( d_t \pmb{u}^n-\dot{\pmb{u}}^n),\eta_h)+\frac{1}{\lambda}(d_t \xi^n-\dot{\xi}^n,\eta_h)-\frac{1}{\lambda}(\vec{\alpha}^{\mathsf{T}} (d_t \vec{p}^n-\dot{\vec{p}}^n),\eta_h)=0.
\end{align}
The combination of \eqref{three_field_form_wp2}, \eqref{coupled_alg2} and \eqref{Stokes_projection} implies that
\begin{align} \label{lemma5.7_eq3}
(\operatorname{div}e_{\pmb{u}}^{h,n},\eta_h)+\frac{1}{\lambda}(e_{\xi}^{I,n}+e_{\xi}^{h,n},\eta_h)-\frac{1}{\lambda}(\vec{\alpha}^{\mathsf{T}}(e_{\vec{p}}^{I,n}+e_{\vec{p}}^{h,n}),\eta_h )=0,\\
\label{lemma5.7_eq4}
   (\operatorname{div}e_{\pmb{u}}^{h,n-1},\eta_h)+\frac{1}{\lambda}(e_{\xi}^{I,n-1}+e_{\xi}^{h,n-1},\eta_h)-\frac{1}{\lambda}(\vec{\alpha}^{\mathsf{T}}(e_{\vec{p}}^{I,n-1}+e_{\vec{p}}^{h,n-1}),\eta_h )=0.
\end{align}
Using \eqref{lemma5.7_eq4} and \eqref{lemma5.7_eq3}, we get
\begin{align} \label{coupled_alg_conv10}
      (\operatorname{div} (d_t e_{\pmb{u}}^{h,n}),\eta_h)+\frac{1}{\lambda}(d_t e_{\xi}^{I,n}+d_t e_{\xi}^{h,n},\eta_h)-\frac{1}{\lambda}(\vec{\alpha}^{\mathsf{T}}(d_t e_{\vec{p}}^{I,n}+d_t e_{\vec{p}}^{h,n}),\eta_h )=0.
\end{align}
Combining \eqref{coupled_alg_conv9} with \eqref{coupled_alg_conv10}, we  obtain
\begin{align} \label{lemma5.7_eq5}
    &(\operatorname{div}(d_t e_{\pmb{u}}^{h,n}),\eta_h)+\frac{1}{\lambda}(d_t e_{\xi}^{h,n},\eta_h)-\frac{1}{\lambda}(\vec{\alpha}^{\mathsf{T}} d_t e_{\vec{p}}^{h,n},\eta_h)=
    (\operatorname{div}(d_t \pmb{u}^n-\dot{\pmb{u}}^n),\eta_h)\\ \notag
    &+\frac{1}{\lambda}(d_t \xi^n-\dot{\xi}^n,\eta_h)
    -\frac{1}{\lambda}(d_t e_{\xi}^{I,n},\eta_h)
    -\frac{1}{\lambda}(\vec{\alpha}^{\mathsf{T}} (d_t \vec{p}^n-\dot{\vec{p}}^n),\eta_h)+\frac{1}{\lambda}(\vec{\alpha}^{\mathsf{T}} d_t e_{\vec{p}}^{I,n},\eta_h).
\end{align}
The above equation \eqref{lemma5.7_eq5} can also be written as
\begin{align} \label{coupled_alg_conv11}
    &(\operatorname{div}(d_t e_{\pmb{u}}^{h,n}),\eta_h)+\frac{1}{\lambda}(d_t e_{\xi}^{h,n},\eta_h)-\frac{1}{\lambda}(\vec{\alpha}^{\mathsf{T}} d_t e_{\vec{p}}^{h,n},\eta_h)=
    (\operatorname{div}(d_t \pmb{u}^n-\dot{\pmb{u}}^n),\eta_h)\\ \notag
    &+\frac{1}{\lambda}(d_t\Pi_h^W \xi^n-\dot{\xi}^n,\eta_h)-\frac{1}{\lambda}(\vec{\alpha}^{\mathsf{T}} (d_t \Pi_h^{M}\vec{p}^n-\dot{\vec{p}}^n),\eta_h).
\end{align}

In \eqref{three_field_form_wp3}, we let $t=t_n$ and $\vec{q}=\vec{q}_h$,  then subtract equation \eqref{coupled_alg3} from it. Combining  with  \eqref{Stokes_projection} and \eqref{ell_projection_op}, we have
\begin{align} \label{coupled_alg_conv12}
    &\left(  \left(  S+\frac{1}{\lambda}\vec{\alpha}\vec{\alpha}^{\mathsf{T}} \right)d_t e_{\vec{p}}^{h,n},\vec{q}_h  \right)-\frac{1}{\lambda}\left( \vec{\alpha}d_t e_{\xi}^{h,n},\vec{q}_h  \right)+\left( K\nabla e_{\vec{p}}^{h,n}, \nabla \vec{q}_h \right) \\ \notag
    &+ \left(  Be_{\vec{p}}^{h,n},\vec{q}_h \right)=\left(  \left(  S+\frac{1}{\lambda}\vec{\alpha}\vec{\alpha}^{\mathsf{T}} \right) \left(d_t \Pi_h^{M}\vec{p}^n-\dot{\vec{p}}^n   \right), \vec{q}_h  \right)-\frac{1}{\lambda}\left( \vec{\alpha}\left( d_t \Pi_h^W \xi^n-\dot{\xi}^n \right)  , \vec{q}_h \right)\\ \notag
&+ \left(     B\left( \Pi_h^{M} \vec{p}^n-\vec{p}^n \right) ,\vec{q}_h\right).
\end{align}
In \eqref{coupled_alg_conv8}, \eqref{coupled_alg_conv11} and \eqref{coupled_alg_conv12}, letting $\pmb{v}_h=d_t e_{\pmb{u}}^{h,n}$, $\eta_h=e_{\xi}^{h,n}$, $\vec{q}_{h}=e_{\vec{p}}^{h,n}$,  then adding them up, one can get
\begin{align*}
 &  2 \mu (\epsilon(e_{\pmb{u}}^{h,n}),\epsilon(d_t e_{\pmb{u}}^{h,n}))+\frac{1}{\lambda}(d_t e_{\xi}^{h,n},e_{\xi}^{h,n})-\frac{1}{\lambda}(\vec{\alpha}^{\mathsf{T}} d_t e_{\vec{p}}^{h,n},e_{\xi}^{h,n})
  +\left(  \left(  S+\frac{1}{\lambda}\vec{\alpha}\vec{\alpha}^{\mathsf{T}} \right)d_t e_{\vec{p}}^{h,n},e_{\vec{p}}^{h,n}  \right)\\ \notag
  & -\frac{1}{\lambda}\left( \vec{\alpha}d_t e_{\xi}^{h,n},e_{\vec{p}}^{h,n}  \right)
  +\left( K\nabla e_{\vec{p}}^{h,n}, \nabla e_{\vec{p}}^{h,n} \right) + \left(  Be_{\vec{p}}^{h,n},e_{\vec{p}}^{h,n} \right)
   =(\operatorname{div}(d_t \pmb{u}^n-\dot{\pmb{u}}^n),e_{\xi}^{h,n})\\ \notag
   &+\frac{1}{\lambda}(d_t\Pi_h^W \xi^n-\dot{\xi}^n,e_{\xi}^{h,n})-\frac{1}{\lambda}(\vec{\alpha}^{\mathsf{T}} (d_t \Pi_h^{M} \vec{p}^n-\dot{\vec{p}}^n),e_{\xi}^{h,n}) -\frac{1}{\lambda}\left( \vec{\alpha}\left( d_t \Pi_h^W \xi^n-\dot{\xi}^n \right)  , e_{\vec{p}}^{h,n} \right)
   \\ \notag
   &+\left(  \left(  S+\frac{1}{\lambda}\vec{\alpha}\vec{\alpha}^{\mathsf{T}} \right) \left(d_t \Pi_h^{M}\vec{p}^n-\dot{\vec{p}}^n   \right), e_{\vec{p}}^{h,n}  \right)
   + \left(     B\left( \Pi_h^{M} \vec{p}^n-\vec{p}^n \right) ,e_{\vec{p}}^{h,n}\right).
\end{align*}
The above formula can also be written as
\begin{align} \label{coupled_alg_conv13}
    & d_t \left( \mu|\!|\epsilon(e_{\pmb{u}}^{h,n})|\!|_{L^2(\Omega)}^2+\frac{1}{2\lambda}|\!|\vec{\alpha} ^{\mathsf{T}} e_{\vec{p}}^{h,n}-e_{\xi}^{h,n}|\!|_{L^2(\Omega)}^2+\frac{1}{2}|\!|\hat{S}e_{\vec{p}}^{h,n}|\!|_{L^2(\Omega)}^2 \right) \\ \notag
   &+|\!|\hat{K}\nabla e_{\vec{p}}^{h,n}|\!|_{L^2(\Omega)}^2+\left(  Be_{\vec{p}}^{h,n},e_{\vec{p}}^{h,n} \right)\\ \notag
&+\Delta t\left(\mu|\!|d_t \epsilon(e_{\pmb{u}}^{h,n})|\!|_{L^2(\Omega)}^2
+\frac{1}{2\lambda}|\!|d_t(\vec{\alpha}^{\mathsf{T}} e_{\vec{p}}^{h,n}-e_{\xi}^{h,n})|\!|_{L^2(\Omega)}^2+\frac{1}{2}|\!|d_t\hat{S} e_{\vec{p}}^{h,n}|\!|_{L^2(\Omega)}^2\right)\\ \notag
&=\left(\operatorname{div}(d_t \pmb{u}^n-\dot{\pmb{u}}^n),e_{\xi}^{h,n}\right)
+\frac{1}{\lambda}\left(d_t \Pi_h^W \xi^n-\dot{\xi}^n,e_{\xi}^{h,n}\right)
-\frac{1}{\lambda}\left(\vec{\alpha}^{\mathsf{T}} (d_t \Pi_h^{M} \vec{p}^n-\dot{\vec{p}}^n),e_{\xi}^{h,n}\right) \\ \notag
&+\left(  \left(  S+\frac{1}{\lambda}\vec{\alpha}\vec{\alpha}^{\mathsf{T}} \right) \left(d_t \Pi_h^{M}\vec{p}^n-\dot{\vec{p}}^n   \right), e_{\vec{p}}^{h,n}  \right)
-\frac{1}{\lambda}\left( \vec{\alpha}\left( d_t \Pi_h^W \xi^n-\dot{\xi}^n \right)  , e_{\vec{p}}^{h,n} \right)\\ \notag
 &+ \left(     B\left( \Pi_h^{M} \vec{p}^n-\vec{p}^n \right) ,e_{\vec{p}}^{h,n}\right).
\end{align}
Applying the summation operator $\Delta t \sum_{n=1}^l$ to both sides, we obtain \eqref{coupled_alg_conv7}.
$\hfill\qedsymbol$

\vspace{0.5cm}
In the rest part of the analysis, we will use the following discrete Gronwall inequality \cite {ahmed2018really}. 

\begin{lemma}
    Let $\tau,\, B$ and $a_j,\, b_j,\, c_j,\, \gamma_j,\, \forall j\geq 1$ be non-negative numbers such that
\begin{align}
    a_n +\sum_{j=1}^n b_j \leq B+\tau \sum_{j=1}^n c_j +\tau \sum_{j=1}^n \gamma_j a_j .
\end{align}
If $\tau \gamma_j \leq 1, j=1,\cdots,n$, there holds
\begin{align}
    a_n +\sum_{j=1}^n b_j \leq C \left( B+\tau \sum_{j=1}^n c_j \right),
\end{align}
where $C=e^{\tau \sum_{j=1}^n \frac{\gamma_j}{1-\tau \gamma_j}}$.
\end{lemma}

\vspace{0.5cm}

\begin{theorem}
    Let $\{(\pmb{u}_h^n,\xi_h^n,\vec{p}_{h}^n)\}$ be the solution of problem \eqref{coupled_alg1}-\eqref{coupled_alg3}, then the following error estimate holds
\begin{align} \label{coupled_alg_conv14}
  & E_h^l+\frac{1}{2}\Delta t\sum_{n=1}^l  |\!|\hat{K}\nabla  e_{\vec{p}}^{h,n}|\!|_{L^2(\Omega)}^2\leq C_1 (\Delta t)^2+C_2 h^4+C_3\Delta t h^4,
\end{align}
where
\begin{align*}
&C_1=C_1 \left(|\!|\ddot{\pmb{u}}|\!|_{L^2(0,t_l;H^1(\Omega))},|\!|\ddot{\xi}|\!|_{L^2(0,t_l;L^2(\Omega))},|\!|\ddot{\vec{p}}|\!|_{L^2(0,t_l;L^2(\Omega))}\right),\\
&C_2=C_2\left(|\!|\dot{\pmb{u}}|\!|_{L^2(0,t_l;H^3(\Omega))},|\!|\dot{\xi}|\!|_{L^2(0,t_l;H^2(\Omega))},|\!|\dot{\vec{p}}|\!|_{L^2(0,t_l;H^2(\Omega))}\right),\\
&C_3=C_3 \left(|\!|\vec{p}|\!|_{L^{\infty}(0,t_l;H^2(\Omega))}\right).
\end{align*}
\end{theorem}

\noindent{\bf Proof}. According to the Cauchy-Schwarz inequality, the  Young inequality with a $ \epsilon_1>0$, and  Taylor's formula, we obtain
\begin{align} \label{coupled_alg_conv15}
    &\Delta t\sum_{n=1}^l \left(\operatorname{div}(\text{d}_t \pmb{u}^n-\dot{\pmb{u}}^n),e_{\xi}^{h,n}\right)\\ \notag
   & =\sum_{n=1}^l \left(\operatorname{div}(\pmb{u}^n-\pmb{u}^{n-1}-\Delta t \dot{\pmb{u}}^n),e_{\xi}^{h,n}\right)\\ \notag
 &\leq \sum_{n=1}^l |\!|\pmb{u}^n-\pmb{u}^{n-1}-\Delta t \dot{\pmb{u}}^n |\!|_{H^1(\Omega)} |\!|e_{\xi}^{h,n}|\!|_{L^2(\Omega)}   \\ \notag
 &\leq \left(\sum_{n=1}^l \frac{(\Delta t)^3}{3}\int_{t_{n-1}}^{t_n} |\!|\ddot{\pmb{u}}|\!|_{H^1(\Omega)}^2 \text{d}s\right)^{\frac{1}{2}} \left(\sum_{n=1}^l |\!|e_{\xi}^{h,n}|\!|_{L^2(\Omega)}^2\right)^{\frac{1}{2}}\\ \notag
 &\leq \left(\frac{(\Delta t)^2}{3}\int_{t_0}^{t_l} |\!|\ddot{\pmb{u}}|\!|_{H^1(\Omega)}^2 \text{d}s\right)^{\frac{1}{2}}\left(\Delta t\sum_{n=1}^l |\!|e_{\xi}^{h,n}|\!|_{L^2(\Omega)}^2\right)^{\frac{1}{2}}\\ \notag
 &\leq \frac{(\Delta t)^2}{6 \epsilon_1}|\!|\ddot{\pmb{u}}|\!|_{L^2(0,t_l;H^1(\Omega))}^2+\frac{\epsilon_1 \Delta t}{2}\sum_{n=1}^l |\!|e_{\xi}^{h,n}|\!|_{L^2(\Omega)}^2.
\end{align}
Similarly, the following term can be bounded by 
\begin{align} \label{coupled_alg_conv16}
  & \Delta t \sum_{n=1}^l \frac{1}{\lambda} (\text{d}_t \Pi_h^W \xi^n-\dot{\xi}^n,e_{\xi}^{h,n})\\ \notag
   &\leq \frac{1}{\lambda^2 \epsilon_2}\left(   Ch^4  |\!|\dot{\xi}|\!|_{L^2(0,t_l;H^2(\Omega))}^2 +Ch^4|\!|\dot{\pmb{u}}|\!|_{L^2(0,t_l;H^3(\Omega))}^2+\frac{(\Delta t)^2}{3}|\!|\ddot{\xi}|\!|_{L^2(0,t_l;L^2(\Omega))}^2      \right) \\ \notag
   &+\frac{\epsilon_2 \Delta t}{2} \sum_{n=1}^l |\!|e_{\xi}^{h,n}|\!|_{L^2(\Omega)}^2.
\end{align}
In consideration of \eqref{projection_err3}, there holds
\begin{align} \label{coupled_alg_conv17}
   &\Delta t \sum_{n=1}^l \frac{1}{\lambda}\left( \vec{\alpha}^{\mathsf{T}} (\text{d}_t \Pi_h^{M} \vec{p}^n-\dot{\vec{p}}^n),e_{\xi}^{h,n} \right) \\ \notag
&\leq \frac{1}{\lambda^2 \epsilon_3} |\!|\vec{\alpha}|\!|^2 \left( Ch^4 |\!|\dot{\vec{p}}|\!|_{L^2(0,t_l;H^2(\Omega))}^2+\frac{(\Delta t)^2}{3}|\!|\ddot{\vec{p}}|\!|_{L^2(0,t_l;L^2(\Omega))}^2 \right)+\frac{\epsilon_3 \Delta t}{2}\sum_{n=1}^l |\!|e_{\xi}^{h,n}|\!|_{L^2(\Omega)}^2.
\end{align}
By use of estimate \eqref{projection_err3} and the Poincar$\rm{\acute{e}}$ inequality, we get
\begin{align} \label{coupled_alg_conv18}
   & \Delta t \sum_{n=1}^l  \left(S\left( \text{d}_t \Pi_h^{M}\vec{p}^n-\dot{\vec{p}}^n\right),e_{\vec{p}}^{h,n}  \right)\\ \notag
 &\leq \frac{|\!|S|\!|^2}{\epsilon_4} \left(Ch^4|\!|\dot{\vec{p}}|\!|_{L^2(0,t_l;H^2(\Omega))}^2+\frac{(\Delta t)^2}{3}|\!|\ddot{\vec{p}}|\!|_{L^2(0,t_l;L^2(\Omega))}^2\right)+\frac{\Delta t\epsilon_4}{2}\sum_{n=1}^l  |\!|\nabla e_{\vec{p}}^{h,n}|\!|_{L^2(\Omega)}^2.
\end{align}
Likewise, based on the Cauchy-Schwarz inequality, we derive
\begin{align} \label{coupled_alg_conv19}
& \Delta t \sum_{n=1}^l \left(  \frac{1}{\lambda}\vec{\alpha}\vec{\alpha}^{\mathsf{T}}(\text{d}_t \Pi_h^{M} \vec{p}^n-\dot{\vec{p}}^n),e_{\vec{p}}^{h,n}\right) \\ \notag
&\leq \frac{|\!|\vec{\alpha}\vec{\alpha}^{\mathsf{T}}|\!|^2}{\lambda^2 \epsilon_5}\left(Ch^4 |\!|\dot{\vec{p}}|\!|_{L^2(0,t_l;H^2(\Omega))}^2+\frac{(\Delta t)^2}{3}|\!|\ddot{\vec{p}}|\!|_{L^2(0,t_l;L^2(\Omega))}^2\right)+\frac{\Delta t \epsilon_5}{2}\sum_{n=1}^l  |\!|\nabla e_{\vec{p}}^{h,n}|\!|_{L^2(\Omega)}^2.\\
\label{coupled_alg_conv20}
& \Delta t \sum_{n=1}^l  \frac{1}{\lambda}(\vec{\alpha}(\text{d}_t \Pi_h^W \xi^n-\dot{\xi}^n),e_{\vec{p}}^{h,n})  \\ \notag
&\leq \frac{|\!|\vec{\alpha}|\!|^2}{\lambda^2 \epsilon_6}\bigg( Ch^4 |\!|\dot{\xi}|\!|_{L^2(0,t_l;H^2(\Omega))}^2
+Ch^4 |\!|\dot{\pmb{u}}|\!|_{L^2(0,t_l;H^3(\Omega))}^2+\frac{(\Delta t)^2}{3}|\!|\ddot{\xi}|\!|_{L^2(0,t_l;L^2(\Omega))}^2\bigg)\\ \notag
&+ \frac{\Delta t\epsilon_6}{2}\sum_{n=1}^l  |\!|\nabla e_{\vec{p}}^{h,n}|\!|_{L^2(\Omega)}^2.\\
 \label{coupled_alg_conv21}
 & \Delta t \sum_{n=1}^l \left( Be_{\vec{p}}^{I,n},e_{\vec{p}}^{h,n}\right)
  \leq \frac{Cl|\!|B|\!|^2h^4}{2\epsilon_7}\Delta t |\!|\vec{p}|\!|_{L^{\infty}(0,t_l;H^2(\Omega))}^2+\frac{\Delta t\epsilon_7}{2}\sum_{n=1}^l  |\!|\nabla e_{\vec{p}}^{h,n}|\!|_{L^2(\Omega)}^2.
\end{align}
By using the discrete inf-sup condition, we see that
\begin{align*}
    |\!|e_{\xi}^{h,n}|\!|_{L^2(\Omega)} \leq \frac{1}{\beta_0^*} \sup_{\pmb{v}_h\in \pmb{V}_h} \frac{(e_{\xi}^{h,n},\operatorname{div}\pmb{v}_h)}{|\!|\pmb{v}_h|\!|_{H^1(\Omega)}}&=\frac{1}{\beta_0^*} \sup_{\pmb{v}_h \in \pmb{V}_h}\frac{2 \mu (\epsilon(e_{\pmb{u}}^{h,n}),\epsilon(\pmb{v}_h))}{|\!|\pmb{v}_h|\!|_{H^1(\Omega)}} \leq \frac{2 \mu}{\beta_0^*}|\!|\epsilon(e_{\pmb{u}}^{h,n})|\!|_{L^2(\Omega)}.
\end{align*}
So we can choose  $\epsilon_i$,\,$i\in\{1,2,3\}$, small enough such that
\begin{align} \label{coupled_alg_conv22}
     \frac{\epsilon_i \Delta t}{2}\sum_{n=1}^l |\!|e_{\xi}^{h,n}|\!|_{L^2(\Omega)}^2 \leq \frac{\Delta t\mu}{3} \sum_{n=1}^l |\!|\epsilon(e_{\pmb{u}}^{h,n})|\!|_{L^2(\Omega)}^2.
\end{align}
Similarly, select $\epsilon_i$,\,$i\in\{4,5,6,7\}$, small enough to guarantee
    \begin{align} \label{coupled_alg_conv23}
        \frac{\Delta t\epsilon_i}{2}\sum_{n=1}^l |\!|\nabla e_{\vec{p}}^{h,n}|\!|_{L^2(\Omega)}^2\leq \frac{\Delta t}{8} \sum_{n=1}^l |\!|\hat{K}\nabla e_{\vec{p}}^{h,n}|\!|_{L^2(\Omega)}^2.
    \end{align}
Combining \eqref{coupled_alg_conv7}, \eqref{coupled_alg_conv15}-\eqref{coupled_alg_conv23} and the discrete Gronwall inequality, we see \eqref{coupled_alg_conv14}.    $\hfill\qedsymbol$

\vspace{0.5cm}

\begin{theorem}
Let $\{(\pmb{u}_h^n,\xi_h^n,\vec{p}_{h}^n)\}$ be the solution of problem \eqref{coupled_alg1}-\eqref{coupled_alg3}, then the following error estimate holds
\begin{align} \label{coupled_alg_conv24}
 \mu|\!|\epsilon(e_{\pmb{u}}^{l})|\!|_{L^2(\Omega)}^2+\frac{1}{2\lambda}|\!|\vec{\alpha} ^{\mathsf{T}} e_{\vec{p}}^{l}|\!|_{L^2(\Omega)}^2+\frac{1}{2}|\!|\hat{S}e_{\vec{p}}^{l}|\!|_{L^2(\Omega)}^2
  &\leq C_1 (\Delta t)^2+C_2 h^4+C_3\Delta t h^4,\\
\label{coupled_alg_conv25}
   \Delta t \sum_{n=1}^l  |\!|\hat{K}\nabla e_{\vec{p}}^n|\!|_{L^2(\Omega)}^2 &\leq C_1 (\Delta t)^2+C_2 h^4+C_3\Delta t h^2,\\
\label{coupled_alg_conv26}
  |\!|e_{\pmb{u}}^l|\!|_{H^1(\Omega)}^2 &\leq C_1(\Delta t)^2+C_2 h^4+C_3 \Delta t h^4,\\
\label{coupled_alg_conv27}
    |\!|e_{\xi}^l|\!|_{L^2(\Omega)}^2 &\leq C_1(\Delta t)^2+C_2h^4+C_3 \Delta t h^4,
\end{align}
where
\begin{align*}
&C_1=C_1 \big(|\!|\ddot{\pmb{u}}|\!|_{L^2(0,t_l;H^1(\Omega))},|\!|\ddot{\xi}|\!|_{L^2(0,t_l;L^2(\Omega))},|\!|\ddot{\vec{p}}|\!|_{L^2(0,t_l;L^2(\Omega))}\big),\\
&C_2=C_2\big(|\!|\pmb{u}|\!|_{L^{\infty}(0,t_l;H^3(\Omega))},|\!|\xi|\!|_{L^{\infty}(0,t_l;H^2(\Omega))},|\!|\vec{p}|\!|_{L^{\infty}(0,t_l;H^2(\Omega))},|\!|\dot{\pmb{u}}|\!|_{L^2(0,t_l;H^3(\Omega))},\\
& |\!|\dot{\xi}|\!|_{L^2(0,t_l;H^2(\Omega))},   |\!|\dot{\vec{p}}|\!|_{L^2(0,t_l;H^2(\Omega))} \big),\\
&C_3=C_3 \big(|\!|\vec{p}|\!|_{L^{\infty}(0,t_l;H^2(\Omega))}\big).
\end{align*}
\end{theorem}
\noindent{\bf Proof}. By using the Cauchy-Schwarz inequality and discrete inf-sup condition, we have
\begin{align} \label{coupled_alg_conv29}
       |\!|\vec{\alpha}^{\mathsf{T}} e_{\vec{p}}^{h,l}|\!|_{L^2(\Omega)}&\leq |\!|\vec{\alpha}^{\mathsf{T}} e_{\vec{p}}^{h,l}-e_{\xi}^{h,l}|\!|_{L^2(\Omega)}+|\!|e_{\xi}^{h,l}|\!|_{L^2(\Omega)}\\ \notag
      &  \leq |\!|\vec{\alpha}^{\mathsf{T}} e_{\vec{p}}^{h,l}-e_{\xi}^{h,l}|\!|_{L^2(\Omega)}+\frac{2 \mu}{\beta_0}|\!|\epsilon(e_{\pmb{u}}^{h,l})|\!|_{L^2(\Omega)}.
\end{align}
    Combining \eqref{coupled_alg_conv14}, \eqref{coupled_alg_conv29} and \eqref{projection_err1}-\eqref{projection_err3}, we obtain \eqref{coupled_alg_conv24}, \eqref{coupled_alg_conv25} and \eqref{coupled_alg_conv27}. According to the first Korn inequality and \eqref{coupled_alg_conv24}, we derive \eqref{coupled_alg_conv26}.
$\hfill\qedsymbol$
\end{appendices}


\bibstyle{plain}
\bibliography{mybib}

\end{document}